\documentclass[a4paper]{amsart}

\usepackage{esint,ifpdf}

\ifpdf
  \usepackage[final,expansion=alltext,protrusion=alltext]{microtype} 
\fi

\usepackage{tikz}
\usetikzlibrary{matrix,positioning,shapes,arrows,decorations,decorations.shapes}
\pgfarrowsdeclarealias{<}{>}{stealth}{stealth}%
\pgfarrowsdeclaredouble{<<}{>>}{stealth}{stealth}%
\usepackage{tikz-cd}

\tikzset{
  commutative diagrams/.cd,
  arrow style=tikz,
  diagrams={>=stealth},
  shift up/.style={
    to path={([yshift=#1]\tikztostart.east) -- ([yshift=#1]\tikztotarget.west) \tikztonodes}},
  mathdouble/.style={-,double equal sign distance}
}

\usepackage[utf8]{inputenc}
\usepackage[british]{babel}

\usepackage{amsmath,amssymb,amsthm,bbm,mathtools,enumitem}

\mathtoolsset{showonlyrefs}

\usepackage{extpfeil} 
\usepackage{stmaryrd} 

\usepackage{tabularx}
\usepackage[alphabetic]{amsrefs}

\newcommand\cf{\emph{cf.}~}
\newcommand\vq{\emph{v.}~}
\newcommand\eg{\emph{e.g.}~}
\newcommand\ie{\emph{i.e.}~}

\newcommand\via{\emph{via}~}

\newcommand\etc{\emph{etc.}}

\newcommand{\C}{\mathbb{C}}
\newcommand{\R}{\mathbb{R}}
\newcommand{\Z}{\mathbb{Z}}
\newcommand{\N}{\mathbb{N}}

\newcommand{\1}{\mathbbm{1}}

\DeclareMathOperator\Gl{Gl}
\newcommand{\gl}{\mathfrak{gl}}

\newcommand{\g}{\mathfrak{g}}
\newcommand{\h}{\mathfrak{h}}
\newcommand{\p}{\mathfrak{p}}
\newcommand{\mfa}{\mathfrak{a}}

\newcommand{\mf}[1]{\mathfrak{#1}}

\DeclareMathOperator\U{U}

\newcommand{\D}{\mathrm{d}}

\newcommand{\id}{\operatorname{id}}

\newcommand{\into}{\hookrightarrow}

\DeclareMathOperator{\diag}{diag}

\DeclareMathOperator{\Ad}{{Ad}}
\DeclareMathOperator{\ad}{{ad}}

\DeclareMathOperator{\End}{End}
\newcommand{\uEnd}{\underline{\End}}
\DeclareMathOperator{\Hom}{Hom}

\DeclareMathOperator{\Aut}{Aut}
\DeclareMathOperator{\Ind}{Ind}

\DeclareMathOperator{\Sets}{Sets}

\newcounter{toProof1}
\newcounter{toProof2}

\newcounter{labelCount}
\newtheorem{theo}{Theorem}[section]
\newenvironment{THM}[1][\arabic{labelCount}]{
  \begin{theo}
    \label{thm:#1}
    \stepcounter{labelCount}
  }{
  \end{theo}
  \setcounter{toProof1}{\value{section}}
  \setcounter{toProof2}{\value{theo}}
}

\newtheorem{lem}[theo]{Lemma}
\newenvironment{LMM}[1][\arabic{labelCount}]{
  \begin{lem}
    \label{lmm:#1}
    \stepcounter{labelCount}
  }{
  \end{lem}
  \setcounter{toProof1}{\value{section}}
  \setcounter{toProof2}{\value{theo}}
}

\newtheorem{coro}[theo]{Corollary}
\newenvironment{CRL}[1][\arabic{labelCount}]{
  \begin{coro}
    \label{crl:#1}
    \stepcounter{labelCount}
  }{
  \end{coro}
  \setcounter{toProof1}{\value{section}}
  \setcounter{toProof2}{\value{theo}}
}

\newtheorem*{prop*}{Proposition}
\newtheorem{prop}[theo]{Proposition}
\newenvironment{PRO}[1][\arabic{labelCount}]{
  \begin{prop}
    \label{pro:#1}
    \stepcounter{labelCount}
  }{
  \end{prop}
  \setcounter{toProof1}{\value{section}}
  \setcounter{toProof2}{\value{theo}}
}

\theoremstyle{definition}
\newtheorem{defin}[theo]{Definition}
\newenvironment{DFN}[1][\arabic{labelCount}]{
  \begin{defin}
    \label{dfn:#1}
    \stepcounter{labelCount}
  }{\end{defin}}

\newtheorem{exa}[theo]{Example}
\newenvironment{EG}[1][\arabic{labelCount}]{
  \begin{exa}
    \label{eg:#1}
    \stepcounter{labelCount}
  }{\end{exa}}

\newtheorem{rema}[theo]{Remark}
\newenvironment{RM}[1][\arabic{labelCount}]{
  \begin{rema}
    \label{rm:#1}
    \stepcounter{labelCount}
  }{\end{rema}}

\newenvironment{PRF}[1][Proof]{\begin{proof}[#1]}{\end{proof}}

\makeatletter
\def\Set@Scallop[#1]#2#3{{#1}\Parens{#2}{#3}}
\newcommand\DeclareScalableOperator[2]{%
  \expandafter\def\csname#1\endcsname{\@ifnextchar[{{#2}\Set@Scallop}{{#2}\Set@Scallop[{}]}}
}
\makeatother

\def\ev{{\bar 0}}
\def\odd{{\bar 1}}
\newcommand{\defi}{\coloneqq}                   

\newcommand{\fa}{for all }

\newcommand\mathfa[1][{}]{\quad\text{\fa{#1} }}

\newcommand{\scth}{such that }
\newcommand{\AND}{and}

\newcommand\mathtxt[1]{\quad\text{{#1}}\quad}

\newcommand{\nd}{\mathtxt\AND}

\newcommand\vphi{\varphi}
\newcommand\vrho{\varrho}

\newcommand\eps{\varepsilon}
\newcommand\nats{\mathbb{N}}
\newcommand\ints{\mathbb{Z}}
\newcommand\rats{\mathbb{Q}}
\newcommand\reals{\mathbb{R}}
\newcommand\cplxs{\mathbb{C}}

\newcommand\vvoid{\varnothing}
\newcommand\sle{\leqslant}
\newcommand\sge{\geqslant}

\makeatletter
\newcommand\Size[7][1]{
  \ifx#20%
  \def\r@l{}\def\r@m{}\def\r@r{}%
  \else%
  \ifx#21%
  \def\r@l{\bigl}\def\r@r{\bigr}\def\r@m{\bigm}%
  \else%
  \ifx#22%
  \def\r@l{\Bigl}\def\r@r{\Bigr}\def\r@m{\Bigm}%
  \else%
  \ifx#23%
  \def\r@l{\biggl}\def\r@r{\biggr}\def\r@m{\biggm}%
  \else
  \ifx#24%
  \def\r@l{\Biggl}\def\r@r{\Biggr}\def\r@m{\Biggm}%
  \fi%
  \fi%
  \fi%
  \fi%
  \fi%
  \ifx#10%
  \def\r@m{}%
  \fi%
  \r@l#3{#4}\r@m#5{#6}\r@r#7%
}%

\makeatother

\newcommand\Set[3]{
  \Size{#1}{\{}{#2}{|}{#3}{\}}%
}%
\newcommand\Dual[3]{
  \Size[0]{#1}{\langle}{#2}{,}{#3}{\rangle}%
}%
\newcommand\Parens[2]{
  \Size[0]{#1}{(}{#2}{}{}{)}
}
\newcommand\Bracks[2]{
  \Size[0]{#1}{[}{#2}{}{}{]}
}
\newcommand\Braces[2]{
  \Size[0]{#1}{\{}{#2}{}{}{\}}
}

\newcommand\Norm[2]{
  \Size[0]{#1}{\lVert}{#2}{}{}{\rVert}
}
\newcommand\Abs[2]{
  \Size[0]{#1}{\lvert}{#2}{}{}{\rvert}
}

\makeatletter

\newif\if@smallmat
\newif\if@none
\newif\if@paren
\newif\if@brack
\newif\if@brace
\newif\if@vline
\newenvironment{Matrix}[2][1]
{\ifx#20%
  \@smallmattrue%
  \else%
  \@smallmatfalse
  \fi%
  \ifx#11%
  \@nonefalse\@parentrue\@brackfalse\@bracefalse\@vlinefalse%
  \else%
  \ifx#12%
  \@nonefalse\@parenfalse\@bracktrue\@bracefalse\@vlinefalse%
  \else%
  \ifx#13%
  \@nonefalse\@parenfalse\@brackfalse\@bracetrue\@vlinefalse%
  \else%
  \ifx#14%
  \@nonefalse\@parenfalse\@brackfalse\@bracefalse\@vlinetrue
  \else%
  \ifx#15%
  \@nonefalse\@parenfalse\@brackfalse\@bracefalse\@vlinefalse%
  \else%
  \@nonetrue\@parenfalse\@brackfalse\@bracefalse\@vlinefalse%
  \fi%
  \fi%
  \fi%
  \fi%
  \fi%
  \if@smallmat%
  \if@none%
  \begin{smallmatrix}%
    \else%
    \if@paren%
    \bigl(\begin{smallmatrix}%
      \else%
      \if@brack%
      \bigl[\begin{smallmatrix}%
        \else%
        \if@brace%
        \bigl\{\begin{smallmatrix}%
          \else%
          \if@vline%
          \bigl\lvert\begin{smallmatrix}%
            \else%
            \bigl\lVert\begin{smallmatrix}%
              \fi%
              \fi%
              \fi%
              \fi%
              \fi%
              \else%
              \if@none%
              \begin{matrix}%
                \else%
                \if@paren%
                \begin{pmatrix}%
                  \else%
                  \if@brack%
                  \begin{bmatrix}%
                    \else%
                    \if@brace%
                    \begin{Bmatrix}%
                      \else%
                      \if@vline%
                      \begin{vmatrix}%
                        \else%
                        \begin{Vmatrix}%
                          \fi%
                          \fi%
                          \fi%
                          \fi%
                          \fi%
                          \fi}%
                        {\if@smallmat%
                          \if@none%
                        \end{smallmatrix}%
                        \else%
                        \if@paren%
                      \end{smallmatrix}\bigr)%
                      \else%
                      \if@brack%
                    \end{smallmatrix}\bigr]%
                    \else%
                    \if@brace%
                  \end{smallmatrix}\bigr\}%
                  \else%
                  \if@vline%
                \end{smallmatrix}\bigr\rvert%
                \else%
              \end{smallmatrix}\bigr\rVert%
              \fi%
              \fi%
              \fi%
              \fi%
              \fi%
              \else%
              \if@none%
            \end{matrix}%
            \else%
            \if@paren%
          \end{pmatrix}%
          \else%
          \if@brack%
        \end{bmatrix}%
        \else%
        \if@brace%
      \end{Bmatrix}%
      \else%
      \if@vline%
    \end{vmatrix}%
    \else%
  \end{Vmatrix}%
  \fi%
  \fi%
  \fi%
  \fi%
  \fi%
  \fi}%

\makeatother

\def\ger{\mathfrak}

\newcommand\CategoryTypeface{\mathbf}
\def\cat{\CategoryTypeface}

\newcommand\SheafTypeface{\mathcal}
\def\sh{\SheafTypeface}

\DeclareMathOperator\str{\mathrm{str}}

\def\DSO{\DeclareScalableOperator}

\DeclareScalableOperator{BER}{\mathrm{Ber}}
\DeclareScalableOperator{IND}{\mathrm{Ind}}
\DeclareScalableOperator{COIND}{\mathrm{Coind}}
\DeclareScalableOperator{DER}{\mathrm{Der}}
\DeclareScalableOperator{Ct}{\mathcal{C}} 
\DeclareScalableOperator{HOM}{\mathrm{Hom}} 
\DeclareScalableOperator{GHOM}{\underline{\mathrm{Hom}}} 
\DeclareScalableOperator{GEND}{\underline{\mathrm{End}}} 
\DeclareScalableOperator{END}{\mathrm{End}}
\DeclareScalableOperator{AUT}{\mathrm{Aut}}
\DeclareScalableOperator{Uenv}{\mathfrak{U}} 
\DeclareScalableOperator{ABer}{\Abs0{\mathrm{Ber}}}

\DSO{ShHom}{\sh Hom} 
\DSO{ShGHom}{\underline{\sh Hom}} 
\DSO{ShEnd}{\sh End} 
\DSO{ShGEnd}{\underline{\sh End}} 
\DSO{ShDer}{\sh Der} 
\DSO{ShGDer}{\underline{\sh Der}} 

\newcommand\fibint[2][\empty]{\!\sideset{#1}{#2}{\fint}} 

\newcommand\Define[1]{\emph{#1}}

\def\sesTheoremn@me#1{\csname#1name\endcsname\ignorespaces}
\def\REF#1#2{\sesTheoremn@me{#1}\ref{#1:#2}}

\newcommand{\pair}[2][0]{ \Parens#1{ \mf{\lowercase{#2}},\uppercase{#2}_{0} } }

\newcommand{\inv}[1]{ { \mf{\lowercase{#1}},\uppercase{#1}_{0} } }

\title{Spherical representations of Lie supergroups}

\author[Alldridge]
{Alexander Alldridge}
\address{Universit\"at zu K\"oln\\
Mathematisches Institut\\
Weyertal 86-90\\
50931 K\"oln}
\email{alldridg@math.uni-koeln.de}

\author[Schmittner]{Sebastian Schmittner}
\address{Universit\"at zu K\"oln\\
Institut f\"ur Theoretische Physik\\
Z\"ulpicher Stra\ss{}e 77\\
50937 K\"oln}
\email{ses@thp.uni-koeln.de}

\begin{document}

\begin{abstract}
  The classical Cartan--Helgason theorem characterises finite-dimen\-sio\-nal \mbox{spherical} representations of reductive Lie groups in terms of their highest weights. We genera\-lise the theorem to the case of a reductive symmetric super\-group pair $(G,K)$ of even type. Along the way, we compute the Harish-Chandra $c$-function of the symmetric superspace $G/K$. By way of an application, we show that in type $A\mathrm{III}|A\mathrm{III}$, all spherical representations are self-dual.
\end{abstract}

\thanks{The first named author was funded by Deutsche Forschungsgemeinschaft (DFG), grant nos.~DFG ZI 513/2-1 and SFB TR/12 ``Symmetries and Universality in Mesoscopic Systems,'' and the Institutional Strategy of the University of Cologne within the German Excellence Initiative. The second named author was partially supported by a grant from the Bonn-Cologne Graduate School of Physics and Astronomy (BCGS), funded by DFG, and a grant from Deutsche Telekom Stiftung}

\keywords{Cartan--Helgason Theorem, Harish-Chandra $c$-function, Lie supergroup, Riemannian symmetric superspace, spherical representation.}

\subjclass[2010]{Primary 17B15, 22E45; Secondary 22E30, 58C50.}

\maketitle


\section*{Introduction}

Let $G$ be a reductive Lie group (with finite centre) and $K$ a maximal compact subgroup. Among the representations of $G$, the spherical ones are those which contain the trivial $K$-representation. If $V$ is finite-dimensional, irreducible, and spherical, then, as the classical Cartan--Helgason Theorem states, the multiplicity equals one, and $V$ admits a quotient map from a spherical principal series representation. Equivalently, the highest weight vector of $V$ is $M$-invariant, where $M=Z_K(\ger a)$. Moreover, the latter condition can be rephrased in terms of the highest weight of $V.$

Spherical representations occur as submodules of the space of functions on the Riemannian symmetric space $G/K$. Hence, replacing $G/K$ by its compact dual symmetric space $U/K$, one obtains as a corollary of the Cartan--Helgason Theorem a complete description of the isotypic summands in the Peter--Weyl decomposition on the space $L^2(U/K)$ of square-integrable functions in terms of their highest weights. 

In this paper, we study sphericity in the setting of \emph{supergroup} representations. Thus, let $(G,K)$ be a symmetric pair of supergroups, which is reductive and of even type. Let $\ger a$ be an even Cartan subspace and $M\defi Z_K(\ger a)$. (For precise definitions, see below.) As our main result, we prove the following theorem (\REF{thm}{main}).

\newtheorem*{theoa*}{Theorem A}
\begin{theoa*}
  Let $V$ be a finite-dimensional $G$-representation of highest weight $\lambda$. If $V^K\neq0$, that is, $V$ is spherical, then $\dim V^K=1|0$ or $0|1$, and the highest weight vector $v_\lambda$ is invariant under the action of $M$. The converse also holds, provided that the highest weight $\lambda$ is high enough. 
\end{theoa*}

Observe that the highest weights of spherical representations are \emph{atypical}, at least if some odd root restricts to zero on $\ger a$ (which is almost always the case). That the highest weight $\lambda$ be \emph{high enough} means that for all odd positive restricted roots $\alpha$, we have $\Dual0{\lambda}\alpha>0$ if $\Dual0\alpha\alpha=0$, and 
\[
  \lambda_\alpha+m_\alpha+2m_{2\alpha}>0\nd\lambda_\alpha+m_\alpha+m_{2\alpha}+1>0
\]  
otherwise, where $\lambda_\alpha\defi\frac{\Dual0\lambda\alpha}{\Dual0\alpha\alpha}$ and $m_\alpha$ is the super-dimension of the restricted root space $\ger g^\alpha$. This condition is forced upon us by the shifted location of the $c$-function zeros, which is induced by the presence of supersymmetries (see below).

As an application of our main result, we investigate the self-duality of spherical representations. Classically, self-duality can be decided by investigating the action of the longest element of the little Weyl group. In particular, in type $A$III, all spherical representations are known to be self-dual. In the corresponding super case, a transitive group action on positive systems is not available. Surprisingly, the self-duality result holds notwithstanding (\REF{crl}{sd}), as we prove as a corollary to \REF{thm}{main}, by a detailed study of even and odd reflections.

\begin{prop*}
  Let $(\ger g,\ger k)$ be the symmetric superpair $(\ger{gl}^{p+q|r+s},\ger{gl}^{p|r}\times\ger{gl}^{q|s})$. Then all finite-dimensional spherical representations are self-dual.
\end{prop*}

Technically, our main result relies heavily on the super-generalisation of the Harish-Chandra $c$-function, \emph{viz} 
\[
  c(\lambda)\defi\int_{\bar N}e^{-(\lambda+\vrho)(H(\bar n))}\,\Abs0{D\bar n},
\]
and the precise location of its zeros. (Here, the Lie superalgebra of $\bar N$ is the sum of the root spaces $\ger g^\alpha$ for negative restricted roots $\alpha$.)

In order to determine the zero locus, we prove a super Gindikin--Karpelevic formula, by establishing a generalised rank reduction procedure, which, instead of the Weyl group, whose action on positive systems is no longer transitive, uses the combinatorial formalism of `odd' and `even' reflections, and works equally well for even and odd, isotropic and anisotropic restricted roots. 

 Explicitly, the formula takes on the following form (\REF{thm}{full-cfn}).

\newtheorem*{theob*}{Theorem B}
\begin{theob*}
  Let $\Re\Dual0\lambda\alpha>0$ for all $\alpha\in\Sigma^+$, $\Dual0\alpha\alpha\neq0$. The integral $c(\lambda)$ converges, and equals
  \[
  c(\lambda)=
  c_0\prod_{\Dual0\alpha\alpha\neq0}2^{-\lambda_\alpha}\frac{\Gamma(\lambda_\alpha)}{\Gamma\Parens1{\tfrac12(\tfrac{m_\alpha}2+1+\lambda_\alpha)}\Gamma\Parens1{\tfrac12(\tfrac{m_\alpha}2+m_{2\alpha}+\lambda_\alpha)}}\prod_{\Dual0\alpha\alpha=0}\Dual0\lambda\alpha^{-\frac{m_\alpha}2}
  \]
  for some non-zero constant $c_0$, independent of $\lambda$. Here, the product extends over all indivisible positive restricted roots.
\end{theob*}

Thus, in directions corresponding to anisotropic roots, $c(\lambda)$ behaves as one might expect from the classical theory, and exhibits simple poles and zeros at integer values of $\lambda_\alpha$, shifted according to the multiplicity of $\alpha$ and its multiples, which may be negative. However, for isotropic roots, it picks up zeros that are not subject to such a shift. The shifted location of the $c$-function zeros is what compels us to impose an extra condition on the highest weight for the sufficiency in Theorem A. 

\bigskip\noindent
In its present form, the classical Cartan--Helgason Theorem was derived by Helgason \cite{helgason-duality} in his study of conical vectors. According to Knapp \cite{knapp-beyond}, the characterisation of sphericity in terms of highest weights was first given (albeit in an incomplete form) by \'E.~Cartan \cite{cartan-1929}, later corrected by Harish-Chandra \cite{hc-spherical1}, and turned into an equivalence by Sugiura \cite{sugiura-1962}. This explains the appellation, which is attributed to G.~Warner.

Schlichtkrull \cite{schlichtkrull} improved Helgason's formulation of the theorem to include non-trivial one-dimensional $K$-types. Our derivation owes much to his exposition. Further generalisations were obtained by Vinberg \cite{vinberg-cone}, Johnson \cite{johnson-helthm}, Kostant \cite{kostant-hel}, and Camporesi \cite{camporesi-rk1}. The theorem's relation to equivariant compactifications was further investigated by Kor\'anyi \cite{koranyi-hel}, who also gave an alternative proof based on the Poisson transform.

This paper is part of an ongoing research project to develop harmonic analysis on Riemannian symmetric superspaces \cite{hhz-sym}. This study is motivated by applications in physics, where symmetric superspaces arise as the target spaces of non-linear supersymmetric $\sigma$-models, \eg in the spectral theory of disordered systems \cite{zirnbauer96} and the study of topological insulators \cite{SRFL09}.

\bigskip\noindent
The paper is organised as follows. In \REF{sec}{pre}, we set up some preliminaries: the relation of supergroups and supergroup pairs, the definition of (possibly infinite-dimensional) smooth representations of supergroup pairs, and the Iwasawa decomposition of reductive symmetric superpairs of even type. These are necessary to formulate our main result, \REF{thm}{main}, in \REF{sec}{main}. In the same section, we also derive the formulation of the main result on the superalgebra level (\ie in terms of highest weights) from the statement at the level of supergroups. In \REF{sec}{mainProof1}, we prove the necessity in \REF{thm}{main}, that is, the $M$-invariance of the highest weight vectors of spherical representations. As part of the proof, we construct possibly infinite-dimensional induced representations of supergroup pairs, and prove Frobenius reciprocity in this setting. 

The converse statement, namely, that the $M$-invariance of the highest weight vector implies that the representation is spherical, at least if the highest weight is high enough, is proved in \REF{sec}{mainProof2}. To that end, we construct, at the supergroup level, a $K$-invariant vector from an $M$-invariant highest weight vector \via a super Eisenstein integral over $K/M$. As in the classical case, the integral is proportional to the $c$-function of $G/K$, although, due to the known pathologies of the Berezin integral, this fact is more difficult to establish than in the even case. We compute $c(\lambda)$ by means of a rank reduction procedure, which applies both to even and odd reflections. This is an extension of the ideas of Gindikin--Karpelevic that avoids the direct use of the Weyl group, whose action is not transitive in the super case. The `rank one' factors at even roots are known from the classical theory. Those occurring at anisotropic odd roots correspond to the rank one symmetric superspaces of even type, which are studied in \cite{ap-cfn}. At isotropic roots, one obtains factors which do not have a geometric counterpart of this type. Nonetheless, their contribution can be evaluated directly.

Other than in the even case, the $c$-functions develop finitely many zeros at weights whose real part is dominant. Hence, depending on the super-dimension of $K/M$, there may be low lying highest weights that are $M$-invariant, but for which we are unable to prove that the corresponding highest weight representations are spherical. 

Finally, in \REF{sec}{selfdual}, we specialise to a particular symmetric superpair where $\mf g=\mf{gl}^{m|n}$. In this case, we show that all spherical highest weight representations are self-dual (\REF{crl}{sd}), by constructing the analogue of the longest little Weyl group element by the use of odd reflections, and study its action on highest weights. Remarkably, for a certain choice of positive root system, the
action is the same as for the longest little Weyl group element of $\mf{gl}^{m+n}$. 

\section{Preliminaries}
\label{sec:pre}

\subsection{Supergroup pairs}

In what follows, we will freely use the theory of supermanifolds, for which one may consult any of the Refs.~\cites{carmeli-caston-fioresi,deligne-morgan,leites,Man88}. 

Let us briefly fix our notation. We will work in the category of \emph{cs} manifolds of J.~Bernstein, which is a full subcategory of the category of $\cplxs$-superspaces. By definition, a $\cplxs$-superspace is a pair $X=(X_0,\sh O_X)$ where $X_0$ is a topological space and $\sh O_X$ is a sheaf of unital supercommutative superalgebras over $\cplxs$, whose stalks $\sh O_{X,x}$ are local rings with maximal ideal $\ger m_{X,x}$. A morphism $f:X\to Y$ is a pair $(f_0,f^\sharp)$ consisting of a continuous map $f_0:X_0\to Y_0$ and a sheaf map $f^\sharp:f_0^{-1}\sh O_Y\to O_X$, which is local in the sense that $f^\sharp(\ger m_{Y,f_0(x)})\subseteq\ger m_{X,x}$ for any $x$. Global sections $f\in\Gamma(\sh O_X)$ of $\sh O_X$ are called \Define{superfunctions}. Due to the locality, the \Define{value} $f(x)\defi f+\ger m_{X,x}\in\sh O_{X,x}/\ger m_{X,x}$ (usually $=\cplxs$) is defined for any $x$.  \Define{Open subspaces} of a $\cplxs$-superspace $X$ are given by $(U,\sh O_X|_U)$, for any open subset $U\subseteq X_0$.

Whenever $V=V_\ev\oplus V_\odd$ is a \emph{real} super-vector space with a complex structure on $V_\odd$, we define $\sh O_V\defi\sh C^\infty_{V_\ev}\otimes_\cplxs\bigwedge_\cplxs(V_\odd)^*$, where $\sh C^\infty$ denotes the sheaf of complex-valued smooth functions. The space $(V_\ev,\sh O_V)$ is called the \Define{\emph{cs} affine superspace} associated with $V$. By abuse of notation, we also denote it by $V$.

Consider now a superspace $X$ whose underlying topological space $X_0$ is Hausdorff and which admits a cover by open subspaces that are isomorphic to open subspaces of some \emph{cs} affine superspace $V$. Such a $\cplxs$-superspace is called a \Define{\emph{cs} manifold}. (This nomenclature is due to J.~Bernstein \cite{deligne-morgan}.) The full subcategory of $\cplxs$-superspaces, whose objects are the \emph{cs} manifolds, admits finite products. 

Group objects in the category of \emph{cs} manifolds are called \emph{cs} Lie supergroups. They can also be understood in terms of the following definition. 

\begin{DFN}
  \label{DEF:SGP}
  A \emph{\emph{cs} supergroup pair} $(\ger g,G_0)$, is given by a
  real Lie group $G_0$ and a complex Lie superalgebra
  $\mf{g}$, with a smooth action $\Ad=\Ad_G : G_0\to
  \Aut(\mf{g})$ by (even) Lie superalgebra
  automorphisms. We demand that the Lie algebra $\ger g_{\ev,\reals}$ of $G_0$ is a
  real form of $\g_\ev$, that $\Ad$ extends the adjoint action of $G_0$ on $\mf{g}_\ev$, and that $\D \Ad$ is given by the bracket of $\ger g$, restricted to $\mf{g}_{\ev,\reals}\times\mf{g}$.

  A \emph{morphism} $(\mf{g},G_0) \to (\mf{h},H_0)$ of \emph{cs} supergroup pairs is a pair $(\D \phi,\phi_0)$ consisting of a morphism of Lie groups $\phi_0 : G_0 \to H_0$ and a morphism $d\phi:\ger g\to\ger h$ of Lie superalgebras, such that $\D \phi$ extends the differential $\D \phi_0$ of $\phi_0$ and intertwines the $\Ad$ actions, \ie
  \[
  \forall g \in G_0 \;\forall x \in \mf{g} : \;\D \phi(\Ad_{G}(g) x) = \Ad_{H}(\phi(g)) \D \phi(x).
  \]

  A \emph{cs} supergroup pair $(\mf{h},H_0)$ is called \emph{subpair} of $(G_0,
  \mf{g})$ if $H_0$ is a closed subgroup of $G_0$, $\mf{h}$ is a Lie
  subsuperalgebra of $\mf{g}$, and $\forall h \in H_0 : \Ad_{G_0}(h)
  |_{\mf{h}} = \Ad_{H_0}(h)$.
\end{DFN}

The following fact is well-known in the case of real supermanifolds (which arise by replacing $\cplxs$ by $\reals$ in the above definitions), see Ref.~\cite{carmeli-caston-fioresi}. The case of \emph{cs} manifolds is virtually identical.

\begin{LMM}
  \label{LMM:LSGrpPairEquivalence}
  The categories of \emph{cs} supergroup pairs and \emph{cs} Lie supergroups are equivalent.
\end{LMM}

The point of view of supergroup pairs will be useful when considering actions of supergroups on possibly infinite-dimensional vector spaces. On the other hand, for the theory of Berezin integration, it will be much more convenient to work on the level of supergroups. 


\subsection{Smooth maps of locally convex spaces}

In this subsection, we give a brief review of basic facts on smooth maps defined on open subsets of locally convex vector spaces. This will allow us to study and construct possibly infinite-dimensional representations of supergroup pairs. 

\begin{DFN}
  Let $V$ be a $\cplxs$-vector space endowed with a Hausdorff topology. If its 
  topology is generated by a set of seminorms \cite{schaefer}, then we will call $V$ a \emph{locally convex
    vector space}. We will \emph{always} impose the Hausdorff condition.

  A locally convex vector space $V$ which is endowed with a $\Z_2$ grading $V=V_\ev\oplus V_\odd$ by closed subspaces is called a \emph{locally convex super-vector space}. 
\end{DFN}

\begin{DFN}
  Let $V$, $W$ be locally convex vector spaces, $U\subseteq V$ an open set, and let $f:U\to W$ be a map. The \emph{directional derivative} of $f$ at $v \in U$ in the
  direction of $x \in V$ is defined by
  \[
  \partial_x |_v f \defi  (\partial_x f)(v) \defi  \partial_t f(v + t x) \defi
  \lim_{t \to 0} \frac{f(v+t x) - f(v)}{t}
  \]
  whenever this exists. This defines the differential of $f$
  \begin{align}
    \D f : U \times V\to W:(x,v)\mapsto \partial_x |_v f.
  \end{align}
  We call $f$ \emph{continuously differentiable} if the map $\D f$ exists and is continuous. The set of all such maps $f$ is denoted $\Ct[^1]0{U,W}$.

  This definition is iterated as usual to define $f\in\Ct[^k]0{U,W}$ if $\D^k f:U\times V^k\to W$ exists and is continuous. The map $f$ is called \emph{smooth} if $f \in \Ct[^k]0{U,W}$ \fa $k\in\nats$. The set of all such maps is denoted by $\Ct[^\infty]0{U,W}$.

  If $M$ is a smooth manifold and $f:M\to W$ is a map, then $f$ is called \emph{smooth} if $f\circ\vphi\in\Ct[^\infty]0{U,W}$ for each local chart $\vphi:U\to M$. The set of all smooth maps $M\to W$ is denoted by $\Ct[^\infty]0{M,W}$.
\end{DFN}

The following fact is well-known. The proof is given for the reader's convenience.

\begin{LMM}[curry]
  Let $U_1\subseteq\reals^{n_1}$ and $U_2\subseteq\reals^{n_2}$
  be open sets and $W$ a locally convex vector space, whose topology is generated by the family $(\Norm0\cdot_j^W)_{j\in J}$ of seminorms. There is a bijection
  \[
  \Ct[^\infty]0{U_1,\Ct[^\infty]0{U_2,W}}\to\Ct[^\infty]0{U_1\times U_2,W}.
  \]

  Here, for $U\subseteq\reals^n$ open, $\Ct[^\infty]0{U,W}$ is given the locally convex topology generated by the seminorms
  \[
  \Norm0h_{j,x_1,\dotsc,x_k,K}\defi \sup_{x\in K}\Norm0{(\partial_{x_k}\dotsm\partial_{x_1}h)(x)}_j^W
  \]
  for $K\subseteq U$ compact, $j\in J$, and $x_1,\dotsc,x_k\in\reals^n$.
\end{LMM}

\begin{PRF}
  For any open $V\subseteq\reals^m$, let $\Ct0{V,W}$ be the set of continuous maps $h:V\to W$, endowed with the compact-open topology. This topology is generated by the seminorms
  \[
  \Norm0h_{j,K}\defi \sup_{x\in K}\Norm0{h(x)}_j^W,
  \]
  for compact $K\subseteq V$ and $j\in J$. By the definition of the topologies, there is for any open $U\subseteq\reals^n$ a topological embedding
  \begin{equation}\label{eq:co-emb}
    \Ct[^\infty]0{U,W}\to\prod_{k=0}^\infty\Ct0{U\times(\reals^n)^k,W}:h\mapsto(d^kh).
  \end{equation}

  Since open sets $V_1\subseteq\reals^{m_1}$, $V_2\subseteq\reals^{m_2}$ are first-countable, there is a bijection
  \begin{equation}\label{eq:curry}
    \Ct0{V_1,\Ct0{V_2,W}}\to\Ct0{V_1\times V_2,W}:h\mapsto(x\mapsto(y\mapsto h(x)(y))),
  \end{equation}
  as is well-known.

  Consider the evaluation map $e:\Ct[^\infty]0{U_2,W}\times U_2\to W:(f,v)\mapsto f(v)$. Inductively, one shows that it possesses an $n$-th derivative
  \[
  d^ne(f,x)(h_1,v_1,\dotsc)=d^nf(x)(v_1,\dotsc,v_n)+\sum\nolimits_id^{n-1}h_i(v_1,\dotsc,\widehat{v_i},\dotsc,v_n)
  \]
  which is continuous by Equation~\eqref{eq:curry}. Thus, if $f\in\Ct[^\infty]0{U_1,\Ct[^\infty]0{U_2,W}}$, then $g:U_1\times U_2\to W$, defined by $g(x,y)=f(x)(y)$, is smooth, since $g=e\circ(f\times\id)$.

  Conversely, let $g\in\Ct[^\infty]0{U_1\times U_2,W}$, and define $f(x)(y)=g(x,y)$. Then we have $f:U_1\to\Ct[^\infty]0{U_2,W}$. Moreover, for $x\in U_1$, $y\in U_2$, we have the equality $d^n(f(x))(y)=d^n_2g(x,y)$ where $d_j$ denotes derivatives with respect to the $j$th argument. Since $d_2g$ is continuous, Equation~\eqref{eq:curry} implies that
  \[
  x\mapsto d^n(f(x)):U_1\to\Ct0{U_2\times (\reals^{n_2})^n,W}
  \]
  is continuous. By Equation~\eqref{eq:co-emb}, $f:U\to\Ct[^\infty]0{V,W}$ is continuous.

  Inductively, its $k$-th derivative is given by
  \[
  d^kf(x)(v_1,\dotsc,v_k)(y)=d_1^kg(x,y)(v_1,\dotsc,v_k),
  \]
  and this is continuous as a map $U_1\times (\reals^{n_1})^k\to\Ct[^\infty]0{U_2,W}$, by a similar argument. This proves the claim.
\end{PRF}

\subsection{Representations of supergroup pairs}\label{sec:spair-rep}

In this subsection, we define smooth representations of supergroup pairs on locally convex super-vector spaces. In the finite-dimensional case, this definition coincides with the usual notion for supergroups. The more general setting allows us to study induced representations of supergroups. 

\begin{DFN}
  \label{DEF:LCSreps}
  Let $(\mf g, G_0)$ be a \emph{cs} supergroup pair and $V$ a locally convex super-vector space. Denote by $\ger{gl}(V)$ the Lie superalgebra of all linear endomorphisms of $V$ with the supercommutator bracket, and by $\Gl(V)_0$ the group of all invertible \emph{even} linear endomorphisms of $V$.

  Given a group homomorphism $\pi_0:G_0 \to \Gl(V)_0$ and a morphism $\pi_{\mf{g}} : \mf{g} \to \mf{gl}(V)$ of Lie superalgebras, the pair $\pi = (\pi_{0},\pi_{\mf{g}})$ is called a \emph{smooth representation} of $(\mf g,G_0)$ if the following holds:
  \begin{enumerate}[wide]
  \item The $G_0$ action is continuous as a map
    \begin{align}
      G_0 \times V\to V:(g,v) \mapsto \pi_0(g)v.
    \end{align}
  \item All vectors $v\in V$ are smooth, \ie the following maps are smooth,
    \begin{align}
      G_0 \to V:g  \mapsto \pi_0(g)v.
    \end{align}
  \item The $\mf{g}$ action is continuous as a map
    \begin{align}
      \mf{g} \times V \to V:(x,v) \mapsto \pi_{\mf{g}}(x) v.
    \end{align}
  \item The action $\pi_\mf g$ extends the differential $d\pi_0$ of the Lie group action,  \ie 
    \[
    \forall x \in \mf{g}_{0,\R},v \in V: \partial_t \big|_0 \pi_0(e^{t x}) v = \pi_{\mf{g}}(x)v.
    \]
  \item The action $\pi_\mf g$ is equivariant w.r.t.~the adjoint action of the pair, \ie
    \[
    \forall x \in \mf{g},g \in G_0 : \pi_{\mf{g}}(\Ad(g) x) = \pi_0(g)\pi_{\mf{g}}(x)\pi_0(g)^{-1}
    \]
  \end{enumerate}

\end{DFN}

\begin{RM}
  Some results concerning ordinary representations on locally convex vector spaces apply:
  \begin{enumerate}[wide]
  \item Conditions (1) and (2) hold simultaneously if and only if $\pi_0$ defines a smooth action of $G_0$ \cite{neeb10}*{Theorem 4.4}.
  \item Since $(x,v) \mapsto \pi_{\mf{g}}(x) v$ is bilinear, Condition~(3) holds if and only if $\pi_\ger g$ is a smooth $\mf{g}$ action.
  \item Condition (1), (2) and (4) imply together that $\pi_{\mf{g}}|_{\ger g_{\ev,\R}}$ is continuous \cite{neeb10}*{Lemma 4.2}. Hence, Condition (3) mainly concerns the odd part $\mf{g}_\odd$.
  \item Since $\mf g$ is finite-dimensional, Condition (3) is actually equivalent to the following: For any $x\in\mf g$, the operator $\pi_\mf g(x)$ is continuous on $V$.

    Indeed, let $(x_i)$ be a basis of $\mf g$. The coefficients $c_i(x)$ of the expression of $x=\sum_ic_i(x)x_i$ in the given basis depend continously on $x$. Hence, the quantity
    \[
    \pi_\mf g(x)v=\textstyle\sum_ic_i(x)\pi_\mf g(x_i)v
    \]
    depends continuously on $(x,v)$ if all of the operators $\pi_\ger g(x_i)$ are continuous.
  \end{enumerate}
\end{RM}

\begin{DFN}[dualRep]\label{DEF:RepMorphs}
  Let $(\ger g,G_0)$ be a \emph{cs} supergroup pair and $(\pi_0,\pi_\mf g)$, $(\rho_0,\rho_\mf g)$ be smooth representations on $V$ and $W$, respectively.

  \begin{enumerate}[wide]
  \item An even continuous linear map $f:V\to W$ is called a \emph{morphism of $\pair g$ representations} if
    \begin{align}
      f \circ \pi_0(g) = \rho_0(g) \circ f \text{ and }
      f \circ \pi_{\mf{g}}(x) =
      \rho_{\mf{g}}(x) \circ f
    \end{align}
    for all $g \in G_0$ and $x \in \mf g$. The set of all these $f$ is denoted by $ \Hom_{\mf g, G_0}(V, W)$.

  \item The representation $(\pi_0, \pi_{\g})$ is called \emph{irreducible} if $V$ does not contain any non-trivial $\pi_{\g}(\g)$ stable closed graded subspaces.

  \item Let $V$ be finite-dimensional. The \emph{dual representation} $(\pi_0^*, \pi_{\g}^*)$ on $V^*$ is
    \begin{align}
      \pi_0^*(g) (\mu) v  \defi  \mu(\pi_0(g^{-1}) v)\text{ and }
      \pi_{\g}^*(x) (\mu) v  \defi -\mu(\pi_g(x) v)
    \end{align}
    for $g \in G_0$, $x \in \g$, $\mu \in V^*$ and $v \in V$. Note that this may not be continuous. 
  \item The subspace of \emph{$\pair g$-invariants} is
    \[
    V^{\inv g} \defi  \Set1{v \in V}{\pi_{\mf{g}}(\mf{g}) v =0 \text{ and } \pi_{0}(G_0) v = v}.
    \]
  \end{enumerate}
\end{DFN}

In the following lemma, recall the formalism of points from Appendix \ref{app:point}.

\begin{LMM}
  Let $V$ be a finite-dimensional smooth $\pair g$-representation and $G$ the
  Lie supergroup corresponding to $\pair g$. Let $\pi$ denote 
  the corresponding $G$-representation on $V$. 
  Define, \fa \emph{cs} manifolds $S$, the set 
  \[
  V^G(S) \defi\Set1{v \in_S V}{\forall g \in_S G : \pi(g)v =
  v}.
  \]
  Then the functor $V^G$ is represented by the super-vector space $V^\inv G$.
\end{LMM}

\begin{PRF}
  This follows from the equivalence of categories of Lie supergroups and supergroup pairs \cite{carmeli-caston-fioresi}.
\end{PRF}

\subsection{Iwasawa decomposition}
\label{sec:Iwasawa}

In order to apply the methods of harmonic analysis to the study of spherical representations, we will use certain super versions of the Iwasawa decomposition. In this subsection, we recall the relevant facts from \cite{a-hchom}, and extend these slightly to include cases such als $\ger{gl}^{n|n}$.

\begin{DFN}
  Let $(\ger g,\theta)$ be a pair consisting of a Lie superalgebra $\ger g$ and an involution $\theta:\ger g\to\ger g$. Then $(\ger g,\theta)$ is called a \emph{symmetric superpair}. We will denote $\ger k\defi\ker(\theta-\id)$ and $\ger p\defi \ker(\theta+\id)$.

  A symmetric superpair is called \emph{reductive} if $\ger g$ is a semi-simple $\ger g_\ev$-module, $\ger z(\ger g)\subseteq\ger g_\ev$, and there exists a $\ger g$- and $\theta$-invariant non-degenerate even supersymmetric form $b$ on $\ger g$. It is called of \emph{even type} if there is an Abelian subalgebra $\ger a$ of $\ger g$, contained in $\ger p_\ev$ and consisting entirely of semi-simple elements of $\ger g_\ev$, such that $\ger p=[\ger k,\ger a]$.

  A \emph{\emph{cs} form} of $(\ger g,\theta)$ is a $\theta$-invariant real form $\ger g_{\ev,\reals}$ of $\ger g$ which is $b$-non-degenerate for some choice of $\ger g$- and $\theta$-invariant form $b$. We write $\ger k_{\ev,\reals}\defi \ger g_{\ev,\reals}\cap\ger k$ and $\ger p_{\ev,\reals}\defi \ger g_{\ev,\reals}\cap\ger p$.

  Given a \emph{cs} form, a \emph{real even Cartan subspace} is a subspace $\ger a_\reals\subseteq\ger p_{\ev,\reals}$ whose complexification $\ger a$ is an even Cartan subspace of $(\ger g,\theta)$. The form induced by $b$ on $\ger a^*$ will be denoted by $\Dual0\cdot\cdot$.
\end{DFN}

\begin{DFN}
  Let $\ger l$ be a real Lie algebra. Recall \cite{borel-rss}*{Lemma 4.1, Definition 4.2} that $\ger l$ is called \emph{compact} if the following equivalent conditions are fulfilled: the set $\ad\ger l\subseteq\END0{\ger l}$ consists of semi-simple endomorphisms with imaginary spectra; and $\ger l$ is the Lie algebra of a compact real Lie group. More generally, if $\vrho$ is a linear representation on a finite-dimensional real vector space $V$, then $\ger l$ is called \emph{$\vrho$-compact} if $\vrho(\ger l)$ generates a compact analytic subgroup of $\Gl(V)$.

  Denoting by $\ad_\ger g$ the adjoint action of $\ger g_\ev\subseteq\ger g$ on $\ger g$, a \emph{cs} form $\ger g_{\ev,\reals}$ will be called \emph{non-compact} if $\ger u_\ev\defi\ger k_{\ev,\reals}\oplus i\ger p_{\ev,\reals}$ is an $\ad_\ger g$-compact real form of $\ger g_\ev$; here, $\ad_\ger g$ denotes the adjoint action of $\ger g_\ev$ on $\ger g$. The condition means that if the symmetric pair $(\ger g_{\ev,\reals},\theta)$ is the infinitesimal pair of a symmetric pair of Lie groups $(G_0,\theta)$, then the associated symmetric space $G_0/K_0$ is Riemannian of non-compact type.
\end{DFN}

\begin{LMM}[nc-form]
  Let $(\ger g,\theta)$ be a reductive symmetric superpair.

  \begin{enumerate}
  \item If $\ger g=\ger z(\ger g)\oplus[\ger g,\ger g]$, then $(\ger g,\theta)$ admits a non-compact \emph{cs} form. In this case, we say that $(\ger g,\theta)$ is \emph{strongly reductive}.
  \item If $(\ger g,\theta)$ is of even type, then it possesses for every non-compact \emph{cs} form a real even Cartan subspace.
  \end{enumerate}
\end{LMM}

\begin{PRF}
  This is the content of \cite{a-hchom}*{Lemma 1.5}
\end{PRF}

\begin{EG}[glnn]
  If $\ger g=\ger{gl}^{m|n}(\cplxs)$ and $m=p+r$, $n=q+s$, then we may consider $\theta(x)\defi sxs^{-1}$ where $s=s^{-1}\defi\diag(1_p,-1_r|1_q,-1_s)$. Then $(\ger g,\theta)$ is a reductive symmetric superpair. It is of even type if and only if $(p-q)(r-s)\sge0$ \cite{ahz-chevalley}*{4.2}. Moreover, $(\ger g,\theta)$ is strongly reductive if and only if $p+r\neq q+s$. However, $(\ger g,\theta)$ always admits a non-compact \emph{cs} form, namely $\ger g_{\ev,\reals}\defi\ger u(p,r)\oplus\ger u(q,s)$.
\end{EG}

\begin{DFN}
  Let $(\ger g,\theta)$ be a symmetric superpair. A triple $(\ger g,G_0,\theta)$ where $(\ger g,G_0)$ is a \emph{cs} supergroup pair is called a \emph{global \emph{cs} form} of $(\ger g,\theta)$ if the Lie algebra $\ger g_{\ev,\reals}$ of $G_0$ is a \emph{cs} form of $(\ger g,\theta)$, and if $\theta$ is an involutive automorphism of $G_0$ (denoted by the same letter as the given involution on $\ger g$) whose differential is the restriction of $\theta$ to $\ger g_{\ev,\reals}$, \scth
  \[
  \Ad(\theta(g))=\theta\circ\Ad(g)\circ\theta\in\End(\ger g)\mathfa g\in G_0.
  \]

  A global \emph{cs} form $(\ger g,G_0,\theta)$ of $(\ger g,\theta)$ is called \emph{non-compact} if $\ger g_{\ev,\reals}$ is a non-compact \emph{cs} form of $(\ger g,\theta)$ and if $\Ad_\ger g(K_0)\subseteq\uEnd(\ger g)$ is compact. Here, $K_0$ denotes the analytic subgroup of $G_0$ generated by $\ger k_{\ev,\reals}$, and $\Ad_\ger g$ denotes the adjoint representation of $G_0$ on the Lie superalgebra $\ger g$.
\end{DFN}
%
%

\begin{DFN}
  Let $(\ger g,\theta)$ be a reductive symmetric superpair of even type, $\ger g_{\ev,\reals}$ a non-compact \emph{cs} form, and $\ger a_\reals$ a real even Cartan subspace. Set $\ger a\defi\ger a_\reals\otimes\cplxs$.

  We have
  \begin{equation}\label{eq:restrrootdecomp}
    \ger g=\ger m\oplus\ger a\oplus\bigoplus\nolimits_{\lambda\in\Sigma}\ger g^\lambda
  \end{equation}
  where $\ger m\defi\ger z_{\ger k}(\ger a)$ is the centraliser of $\ger a$ in $\ger k$, and for $\lambda\in\ger a^*$,
  \[
  \ger g^\lambda\defi\Set1{x\in\ger g}{\forall h\in\ger a\,:\,[h,x]=\lambda(h)x}\nd\Sigma\defi\Set1{\lambda\in\ger a^*\setminus0}{\ger g^\lambda\neq0}.
  \]
  We also define $\ger g_j^\lambda\defi\ger g_j\cap\ger g^\lambda$ and $\Sigma_j\defi\Set1{\lambda\in\ger a^*\setminus0}{\ger g_j^\lambda\neq0}$. The elements of $\Sigma$ are called \emph{restricted roots}, with those of $\Sigma_\ev$ and $\Sigma_\odd$ being called \emph{even} and \emph{odd}, respectively. A restricted root $\alpha$ is called \Define{indivisible} if $\frac\alpha2\notin\Sigma$.

  We have $\Sigma=\Sigma_\ev\cup\Sigma_\odd$, but the union may not be disjoint. Occasionally, we will write $\Sigma(\ger g:\ger a)=\Sigma$ and $\Sigma(\ger g_j:\ger a)=\Sigma_j$. The even restricted roots $\lambda\in\Sigma_\ev$ are real on $\ger a_\reals$. Let $\ger g_{\ev,\reals}^\lambda\defi\ger g_{\ev,\reals}\cap\ger g_\ev^\lambda$ \fa $\lambda\in\Sigma_\ev$ and $\ger m_{\ev,\reals}\defi\ger z_{\ger k_{\ev,\reals}}(\ger a_\reals)$.

  Let $\Sigma^+\subseteq\Sigma$ be a \emph{positive system}, \ie a subset \scth $\Sigma=\Sigma^+\,\dot\cup\,{-\Sigma^+}$ and $\Sigma\cap(\Sigma^++\Sigma^+)\subseteq\Sigma^+$. Let $\Sigma_j^+\defi\Sigma_j\cap\Sigma^+$. Then $\Sigma_\ev^+$ is a positive system of the root system $\Sigma_\ev$. Set
  \[
  \ger n\defi\bigoplus\nolimits_{\lambda\in\Sigma^+}\ger g^\lambda\nd\ger n_j\defi\ger g_j\cap\ger n.
  \]
  By the assumptions on $\Sigma^+$, $\ger n=\ger n_\ev\oplus\ger n_\odd$ is an $\ger a$-invariant subsuperalgebra. Moreover, $\ger n_{\ev,\reals}\defi\ger g_{\ev,\reals}\cap\ger n$ (which is a real form of $\ger n_\ev$) is an $\ger a_\reals$-invariant nilpotent subalgebra of $\ger g_{\ev,\reals}$. Since the roots in $\Sigma_\ev$ are real on $\ger a_\reals$,
  \[
  \ger n_{\ev,\reals}=\bigoplus\nolimits_{\lambda\in\Sigma_\ev^+}\ger g_{0,\reals}^\lambda.
  \]

  We have $\ger g=\ger k\oplus\ger a\oplus\ger n$. We call this the \emph{Iwasawa decomposition} associated with the positive system $\Sigma^+$. Given a non-compact global \emph{cs} form $(\ger g,G_0,\theta)$, we say that it has a \emph{global Iwasawa decomposition} if the multiplication map
  \[
  K_0\times A\times N_0\to G_0
  \]
  is a diffeomorphism, where $K_0$, $A$, and $N_0$, are the analytic subgroups of $G_0$ associated with $\ger k_{\ev,\reals}$, $\ger a_\reals$, and $\ger n_{\ev,\reals}$, respectively.
\end{DFN}

\begin{PRO}[even-iwasawa]
  Let $(\ger g,\theta)$ be a reductive symmetric superpair of even type with non-compact \emph{cs} form $\ger g_{\ev,\reals}$ and real even Cartan subspace $\ger a_\reals$. Let $(G_0,\ger g,\theta)$ be a non-compact global \emph{cs} form with global Iwasawa decomposition, and $G$, $K$, and $N$ be the \emph{cs} Lie supergroups associated with the \emph{cs} supergroup pairs $(\ger g,G_0)$, $(\ger k,K_0)$ and $(\ger n,N_0)$, respectively. Then the multiplication morphism
  \[
  K\times A\times N\to G
  \]
  is an isomorphism of \emph{cs} manifolds.
\end{PRO}

\begin{PRF}
  This follows from the proof of \cite{a-hchom}*{Proposition 1.11}.
\end{PRF}

\begin{LMM}
  Let $(\ger g,\theta)$ be a strongly reductive symmetric superpair of even type. Given a non-compact \emph{cs} form $\ger g_{\ev,\reals}$ and real even Cartan subspace $\ger a_\reals$, there is a non-compact global \emph{cs} form $(\ger g,G_0,\theta)$ with global Iwasawa decomposition, such that $\ger g_{\ev,\reals}$ is the Lie algebra of $G_0$.
\end{LMM}

\begin{PRF}
  This is the content of \cite{a-hchom}*{Proposition 1.10} and the first part of the proof of Proposition 1.11 (\emph{op.cit.}).
\end{PRF}

\begin{EG}[gl-iwasawa]
  In the case of the symmetric superpair $(\ger g,\theta)$ considered in \REF{eg}{glnn}, a non-compact global \emph{cs} form with global Iwasawa decomposition always exists when the superpair has even type. Namely, one may take $G_0\defi\U(p,r)\times\U(q,s)$ regardless whether $(\ger g,\theta)$ is strongly reductive (\ie $p+r\neq q+s$) or not.
\end{EG}

\begin{PRO}[Iwasawa-exist]
  Let $(\ger g,\theta)$ be a symmetric superpair of even type where $\ger g$ is contragredient and admits a non-degenerate $\ger g$- and $\theta$-invariant even supersymmetric form $b$ such that $\ger k$ is $b$-non-degenerate.

  Then a non-compact global \emph{cs} form with global Iwasawa decomposition exists.
\end{PRO}

\begin{PRF}
  Unless the Cartan matrix of $\ger g$ is of type $A$, $\ger g$ is simple, so that $(\ger g,\theta)$ is strongly reductive. Consider the case of a Cartan matrix of type $A$, so that $\ger g=\ger{gl}^{m|n}(\cplxs)$. For an involution of the form $\theta(x)=sxs^{-1}$ where $s=\diag(1_p,-1_q,1_r,-1_s)$, the assertion has been proved in \REF{eg}{gl-iwasawa}. By \cite{serganova-invol}, there are only two more conjugacy classes of involutions to consider.

  Indeed, if $n$ is even, set $\theta(x)\defi-sx^{st}s^{-1}$, where we let $s\defi\diag(1_m,J_n)$ and $J_n\defi\begin{Matrix}00&1\\-1&0\end{Matrix}$. Thus, $\ger k=\ger{osp}^{m|n}$, and we may take $G_0\defi\Gl(m,\reals)\times\Gl(\frac n2,\mathbb H)$. For in this case, $K_0=\mathrm O(m)\times\mathrm{USp}(n)$, and the assertion follows from \REF{pro}{even-iwasawa}. Finally, if $m=n$, the only remaining conjugacy class of involution is represented by $\theta(x)=\Pi(x)$ where $\Pi(x)$ is the matrix of the action of $x$ in the parity reversed standard basis of $\cplxs^{m|m}$. But for this involution, the fixed algebra does not admit an even non-degenerate invariant form.
\end{PRF}

\begin{DFN}
  Let $(\ger g,\theta)$ be a reductive symmetric superpair, $\ger a$ an even Cartan subspace, and $\ger h\subseteq\ger g_\ev$ a $\theta$-stable Cartan subalgebra of $\ger g$ containing $\ger a$. The set of $\ger h$-roots of $\ger g$ is denoted by $\Delta=\Delta(\ger g:\ger h)$, and we write  $\ger g^\gamma_\ger h$ for the $\ger h$-root space for $\gamma\in\Delta$. A positive system $\Delta^+$  of $\Delta$ is called \emph{compatible} (with $\theta$) if it induces a positive system of $\Sigma$, \ie 
  $(\Delta^+|_\ger a)\setminus0$ is a positive system of $\Sigma$. 

  Let $V$ be a finite-dimensional $G$-representation and $v\in V$,
  $v\neq0$, an $\ger h$-weight vector. We say that $v$ is a
  \emph{highest weight vector} if $\ger g^\beta_\ger h v=0$ for every 
  $\beta\in\Delta^+$. The representation $V$ is called a
  \emph{highest weight representation} if it admits a cyclic
  highest weight vector. If $V$ is irreducible, then $V$ is a highest weight representation
  if and only if it is an $\ger h$-weight module, \ie the
  direct sum of $\ger h$-weight spaces. In this case the highest
  weight space $V^\lambda=\cplxs v$ is one dimensional.
\end{DFN}

\begin{DFN}
  If $\Delta^+$ is a positive system of $\Delta=\Delta(\ger g:\ger h)$, then we will denote the corresponding
  simple system by $B(\Delta^+)\defi \Delta^+\setminus\Parens1{\Delta^++\Delta^+}$. Similarly, if $\Sigma^+$ is a positive system of $\Sigma=\Sigma(\ger g:\ger a)$, then we denote by $B(\Sigma^+)\defi\Sigma^+\setminus\Parens1{\Sigma^++\Sigma^+}$ the corresponding simple system.
\end{DFN}

\subsection{Odd reflections}
\label{sec:oddRefl}
In addition to the Weyl group of $\mf g_\ev$, whose action is induced by the adjoint action of $G_0$ on $\ger g$, we will use
so-called \emph{odd reflections}. In this subsection, we collect some of
their properties, as detailed in Ref.~\cite{cheng_wang}.

\begin{DFN}[oddrefl]
  For an (odd) isotropic positive root, $\alpha \in \Delta^+_\odd$,
  $\Dual0\alpha\alpha=0$, we denote by $r_\alpha$ the change
  of positive root system from $\Delta^+$ to
  $r_\alpha(\Delta^+)\defi \{-\alpha\} \cup \Delta^+
  \setminus\{\alpha\}$. This is called an \emph{odd reflection} with
  respect to $\alpha$.
\end{DFN}

\begin{LMM}
  If $\Pi=B(\Delta^+)$ is the simple system for $\Delta^+$, then
  \[
  r_\alpha(\Pi) \defi  \Set1{\beta \in \Pi \setminus\{\alpha\}}{\Dual0\beta\alpha = 0} \cup
  \Set1{\beta + \alpha}{\beta \in \Pi \text{ and } \Dual0\beta\alpha \neq 0} \cup \Braces1{-\alpha} 
  \]
  is the one for $r_\alpha(\Delta^+)$.
\end{LMM}
\begin{PRF}
  See \cite{cheng_wang}*{Section 1.3.6, Lemma 1.26}.
\end{PRF}

\begin{DFN}[oddReflOfHW]
  For an odd reflection $r_{\alpha}$ and any highest weight $\lambda \in \mf{h}^*$ of a finite-dimensional simple $\ger g$-module, we define
  \[
  r_\alpha (\lambda) \defi
  \begin{cases}
    \lambda & \Dual0\lambda\alpha=0,\\
    \lambda - \alpha & \Dual0\lambda\alpha \neq 0.
  \end{cases}
  \]
  For a highest weight representation $V$ with highest weight
  space $V^\lambda$, we define
  \[
  r_\alpha (V^\lambda) \defi
  \begin{cases}
    V^\lambda & \Dual0\lambda\alpha=0,\\
    \mf{g}^{-\alpha}V^\lambda & \Dual0\lambda\alpha \neq 0.
  \end{cases}
  \]
\end{DFN}

\begin{LMM}
  Let $V^\lambda$ be the $\Delta^+$-highest weight space of a finite dimensional
  irreducible representation $V$ and
  $\alpha \in B(\Delta^+)_\odd$ an isotropic odd simple root.
  Then the $r_\alpha(\Delta^+)$-highest weight space of $V$ is $r_\alpha(V^\lambda)=V^{r_\alpha(\lambda)}$.
\end{LMM}

\begin{PRF}
  See \cite{cheng_wang}*{Lemma 1.36}.
\end{PRF}

\begin{DFN}
  We call $R=r_{\alpha_1} \circ \ldots \circ r_{\alpha_n}$ a \Define{chain of simple reflections} with respect to a set $\Pi$ of simple roots if
  \[
  \forall i: \alpha_i \in r_{\alpha_{i+1}}(\dotsm r_{\alpha_{n-1}}(r_{\alpha_n} (\Pi))\dotsm).
  \]
\end{DFN}

\begin{RM}
  Note that \emph{simple} reflections can be iterated. That is, for a chain
  of simple reflections $R$ with respect to $\Pi$, and $V^\lambda$ denoting the
  highest weight space of $V$ with respect to $\Pi$, we have that $R(V^\lambda)$ is the
  highest weight space of $V$ of weight $R(\lambda)$ with respect to $R(\Pi)$.
\end{RM}


\section{Statement of the main theorem}
\label{sec:main}

In the following, we consider a reductive symmetric superpair $(\ger
g,\theta)$ of even type as in \REF{pro}{Iwasawa-exist}. Given a non-compact global
\emph{cs} form $(\ger g,G_0,\theta)$ with global Iwasawa
decomposition, a real even Cartan subspace $\ger a_\reals$, and a
positive system $\Sigma^+\subseteq\Sigma=\Sigma(\ger g:\ger a)$, we
will write $\ger g=\ger k\oplus\ger a\oplus\ger n$ for the
corresponding Lie superalgebra Iwasawa decomposition, and similarly
$G_0=K_0AN_0$ for the Iwasawa decomposition on the group
level. Moreover, we let $\ger m\defi\ger z_\ger k(\ger a)$, $M_0\defi
Z_{K_0}(\ger a)$, and denote by $G$, $K$, $M$, and $N$ the \emph{cs}
Lie supergroups corresponding to the \emph{cs} supergroup pairs $(\ger
g, G_0)$, $\pair k$, $\pair m$, and $\pair n$, respectively. Moreover,
let $\mf{h}\subseteq \mf{g}$ be a $\theta$-invariant Cartan subalgebra
containing $\ger a$.

\begin{DFN}
  Let $V$ be a finite-dimensional representation of $G$. Then $V$ is called \emph{spherical} if $V^K\neq0$.
\end{DFN}

\begin{DFN}[highEnough]
  For any anisotropic $\alpha\in\Sigma$, denote
  \[
    \lambda_\alpha\defi\Dual0{\lambda}{\alpha}/\Dual0{\alpha}{\alpha}\nd m_\alpha\defi\dim\mf g_\ev^\alpha-\dim\ger g^\alpha_\odd.
  \]
 
 We call a highest weight $\lambda$ \Define{high enough} if \fa isotropic positive 
  $\beta\in\Sigma^+$, we have $\Dual0{\lambda}{\beta}>0$, and for all odd anisotropic indivisible  $\alpha\in
  \Sigma^+$, we have
  \[
    \lambda_\alpha+m_\alpha+2m_{2\alpha}>0\nd\lambda_\alpha+m_\alpha+m_{2\alpha}+1>0.
  \]  
  Note that isotropic restricted roots are purely odd, but that anisotropic restricted roots may be purely even, purely odd, or odd and even simultaneously. 
\end{DFN}

The following is our main result. 

\begin{THM}[main]
  Assume that $(\mf g, \theta)$ admits a non-compact global \emph{cs}
  form $(\ger g,G_0,\theta)$ with global Iwasawa decomposition. Let
  $V$ be a finite-dimensional irreducible highest weight
  $G$-representation with highest weight $\lambda$.
\begin{enumerate}[wide]
\item If $V$ is spherical, then 
$V^N\simeq V^K \simeq \C$ as smooth $M$-representations.
\item  If $V^N$ is $M$-invariant and
  $\lambda$ is high enough, then $V$ is spherical.
\end{enumerate}
\end{THM}

\begin{RM}
Notice that there might be spherical representations, \eg the
trivial one, for which the highest weight is not high enough.
\end{RM}

The \emph{proof} of the two parts will be given separately in \REF{sec}{mainProof1}
and \REF{sec}{mainProof2}. We will first give a number of
corollaries. 

\begin{CRL}
A finite-dimensional $G$-representation $V$ whose highest weight is high enough is spherical if and only if it has an $M$-invariant highest weight vector.
\end{CRL}

The $M$-invariance of the highest weight may be characterised algebraically.

\begin{PRO}[conditionsOnHW]
  Let $\lambda \in \mf{h}^*$ such that $(V,\pi)$ is a finite-dimensional highest weight module with highest weight vector
  $v_\lambda$ of weight $\lambda$. Then $v_\lambda \in V^M$ if and only if
  \begin{enumerate}
  \item $\lambda \big |_{\mf{k} \cap \mf{h}} \equiv 0$  \label{ITEM:ConditionsOnHW-1}
  \item $\forall \alpha \in \Sigma_\ev^+ :\lambda_\alpha\in
    \N$ \label{ITEM:ConditionsOnHW-3}
  \end{enumerate}
\end{PRO}

\begin{PRF}
  We have $v_\lambda \in V^M$ if and only if $\pi_{\g}(\mf{m}) v_\lambda=0$ and  $\pi_{0}(M_0) v_\lambda=v_\lambda$. Furthermore,
  \[
  \mf{m} =  \mf{k} \cap \mf{h} \oplus \bigoplus_{\beta \in \Delta_{-}}\ger g^\beta_\ger h
  \]
  where $\Delta_-\subseteq\Delta$ is the set of roots vanishing on $\ger a$. Thus, $\pi_{\g}(\mf{m}) v_\lambda =0 $ implies (\ref{ITEM:ConditionsOnHW-1}) and the following condition:
  \begin{equation}\label{ITEM:ConditionsOnHW-2}
    \pi_{\g}(\ger g^{-\beta}) v_\lambda =0\mathtxt{\fa isotropic} \beta \in \Delta_-\cap\Delta^+_\odd.
  \end{equation}

  Conversely, we show that (\ref{ITEM:ConditionsOnHW-1}) and \eqref{ITEM:ConditionsOnHW-2} imply that $\pi_{\g}(\mf{m}) v_\lambda =0 $. To that end, note that for $\beta \in \Delta^+$, $\pi_{\g}(\g^{\beta}) v_\lambda \subseteq V_{\lambda +
    \beta} = 0$ because $\lambda$ is a highest weight. Let $\beta \in
  \Delta_-\cap\Delta^+$ be even. The element $s_\beta$ in the Weyl group $W=W(\ger g:\ger h)$ of $\g_\ev$ satisfies $\lambda - \beta = s_\beta ( \lambda +
  \beta)$ since $\lambda|_{\mf{k}\cap\mf{h}}=0$, and hence $\Dual0\lambda\beta=0$.
  But since the weights of $V$ are $W$-stable, we have $0=V^{\lambda-\beta}\supseteq\pi_{\g}(\ger g^{-\beta})$. Similarly, if $\beta$ is an anisotropic odd root in $\Delta_-\cap\Delta^+$, then $2\beta$ is an even root in $\Delta_-\cap\Delta^+$, and we have $\lambda-\beta=s_{2\beta}(\lambda+\beta)$, so that $\pi_\ger g(\ger g^{-\beta})v_\lambda=0$ in this case, too. Due to condition \eqref{ITEM:ConditionsOnHW-2}, we find $\pi_\ger g(\ger m)v_\lambda=0$, as claimed.

  Now we show that \eqref{ITEM:ConditionsOnHW-2} is implied by (\ref{ITEM:ConditionsOnHW-1}), so that  (\ref{ITEM:ConditionsOnHW-1}) is equivalent to $\pi_\ger g(\ger m)v_\lambda=0$. First, observe that $\ger k\cap\ger h$ is a Cartan subalgebra for $\ger m$, and that $\Delta_-$ is the root system of $\ger m$. So any root in $\Delta_-\cap\Delta^+$ is the positive linear combination of simple roots of $\Delta_-\cap\Delta^+$. The latter are also simple in $\Delta^+$. Indeed, assume the contrary, that is, we have $\alpha=\beta+\gamma$ for some $\alpha,\beta,\gamma\in\Delta^+$ where $\alpha$ is a $\Delta_-\cap\Delta^+$-simple root; in particular, $\alpha|_\ger a=0$. But then $\beta|_\ger a=-\gamma|_\ger a\neq0$, so that $\pm\beta|_\ger a\in\Sigma^+$. This is a contradiction, because is $\Sigma^+$ is a positive system. Thus, in summary, assuming (\ref{ITEM:ConditionsOnHW-1}), it will suffice to prove \eqref{ITEM:ConditionsOnHW-2} for the isotropic simple roots $\beta\in\Delta_-$.

  So, for any $\alpha\in\Delta^+$, let $h_\alpha\in\ger h$ be defined by $b(h,h_\alpha)=\alpha(h)$. By a standard argument, we may choose for all simple $\alpha$, $\beta$ non-zero $e_\alpha\in\ger g^\alpha$, $f_\beta\in\ger g^{-\beta}$ such that $[e_\alpha,f_\beta]=0$ for $\alpha\neq\beta$ and $[e_\alpha,f_\alpha]=h_\alpha$ otherwise.

  Let $\beta\in\Delta_-$ be an odd isotropic simple root. For any simple root $\alpha$,
  \[
  \pi_\ger g(e_\alpha) \pi_\ger g(f_\beta) v_\lambda = \pi_\ger g([e_\alpha, f_\beta])v_\lambda=\delta_{\alpha\beta}\lambda(h_\beta)v_\lambda=0.
  \]
  If $\pi_\ger g(f_\beta)v_\lambda\in V^{\lambda-\beta}$ were non-zero, then it would be a highest weight vector. But the highest weight of $V$ is unique, contradiction. By virtue of the above arguments, this proves finally that (\ref{ITEM:ConditionsOnHW-1}) implies $\pi_\ger g(\ger m)v_\lambda=0$.

  We now study the invariance under the $M_0$ part. By assumption, $V$ is an $\ger h$-weight module, so it is semisimple as a $\ger g_\ev$-module. Since $G_0$ is connected and the weight space $V^\lambda$ is one-dimensional, the vector $v_\lambda$ is a highest weight vector of an irreducible $G_0$-respresentation $U$ (say) contained in $V$. We assume that $\pi_\ger g(\ger m)v_\lambda=0$. By \cite{helgason84}*{Chapter V, proof of Theorem 4.1}, applied to $U$, we have $\pi_0(M_0)v_\lambda=v_\lambda$ if and only if $\lambda_\alpha\in\nats$ \fa $\alpha\in\Sigma_\ev^+$.
\end{PRF}

\begin{CRL}[conditionsOnHW]
  Let $(\mf g, \mf k)$ be one of the symmetric pairs to which
  \REF{thm}{main} applies. Let $\mf h$ be a Cartan subalgebra and
  $\Delta^+$ a compatible positive root system.
  \begin{enumerate}
  \item
    Let $\lambda \in \h^*$ be the highest weight of the finite-dimensional highest weight representation $V$, $\lambda |_{\mf{h} \cap \mf{k}} \equiv 0$,
    $\forall \alpha\in\Sigma_\ev^+:\lambda_\alpha\in\N$, and $\lambda$ be high enough.
    Then $V$ is spherical.
  \item
    Conversely, let $V$ be a finite-dimensional irreducible
    spherical representation with highest weight $\lambda$.
    Then $\lambda |_{\mf{h} \cap \mf{k}} \equiv 0$ and $\forall \alpha \in \Sigma^+_\ev : \lambda_\alpha \in\N$.
  \end{enumerate}
\end{CRL}

\begin{RM}
  Let $(\g, \theta)$ be a symmetric pair as in \REF{thm}{main} and assume that we have $\Delta_-\cap\Delta_\odd \neq\vvoid$, \ie there is $\alpha \in\Delta_\odd$ \scth $\alpha|_{\mf a}=0$. Then all spherical representations of the pair are atypical.

  Indeed, let $L(\lambda)$ be the simple module of highest weight $\lambda$ and $V^0(\lambda)$ the simple $\ger g_\ev$-module of the same highest weight. Let 
  $K(\lambda)$ denote the module constructed in \cite{kac}*{2.b)} and denoted $\bar V(\lambda)$ there. (For $\ger g$ of type I, this is the Kac module.) Then by \cite{kac}*{Propositions 2.1 and 2.4}, we have inside $K(\lambda)$, that 
  \[
    f_\beta v_\lambda \in f_\beta V^0(\lambda)\setminus\{0\}
  \]
  for any $\beta \in \Delta^+_\odd$, \scth $\beta|_\ger a=0$.

  But then \REF{crl}{conditionsOnHW} implies $\Dual0\lambda\beta=0$, so that 
  \[
    f_\beta v_\lambda \notin V(\lambda)\setminus\{0\},
  \]
  and hence $V(\lambda)\neq K(\lambda)$. That is, $\lambda$ is atypical \cite{kac}*{Theorem 4.1}. 
\end{RM}

\section{Spherical representations have $M$-invariant highest weights}
\label{sec:mainProof1}

In this section, we prove the necessity part of \REF{thm}{main}, stated below as \REF{pro}{mainProof1}.
Just as in the classical case, we will use induced representations, as introduced in \REF{sec}{InducedReps}. Since these are potentially infinite-dimensional, we study them in terms of representations of supergroup pairs, which were defined above in \REF{sec}{spair-rep}. Frobenius reciprocity also holds in this case, as we show in \REF{sec}{frobenius}. Using an embedding of the highest weight representation, constructed in \REF{pro}{embedding}, the assertion of \REF{pro}{mainProof1} follows.

Let $(\mf{g}, G_0)$, $\pair k$, $\pair m$, $\pair a$, and $\pair n$ be as in the statement of \REF{thm}{main}. We let $\pair q$ denote the minimal parabolic subpair, defined by $\mf{q}\defi \mf{m}\oplus \mf{a} \oplus \mf{n}$ and $Q_0\defi M_0A_0N_0$.

\begin{PRO}[mainProof1]
  Let $V$ be a finite-dimensional irreducible smooth $\pair g$-repre\-sen\-tation. Then $V^\inv k$ is at most one-dimensional.

  If $V$ is spherical, then $V^\inv n \simeq V^\inv k$ as smooth $(\mf{m},M_0)$-representations. That is, the highest weight vector is $\pair m$-invariant.
\end{PRO}

We give the \emph{proof} immediately, deferring ancillary definitions and constructions to the sections below. 

\begin{PRF}[Proof of \REF{pro}{mainProof1}]
  First, we embed $V$ into the induced representation
  \begin{align}
    V &\into \Ind^{\mf g,G_0}_{\mf q, Q_0}(V^\inv n )
    \label{thm1re1}
  \end{align}
  \via \REF{pro}{embedding}.
  Then the multiplicity of any irreducible $(\mf k, K_0)$-representation $W$ (say) is
  given by the dimension of
  \begin{align}
    \HOM[_{\mf k, K_0}]1{W, \Ind^{\mf g, G_0}_{\mf q,Q_0}(V^\inv n)}
    &= \Hom_{\mf m, M_0} \left(W, V^\inv N \right)
    \label{thm1re2}
  \end{align}
  by Frobenius reciprocity (\REF{pro}{Frobenius}), which is applicable by \REF{pro}{inducedQGMK}.

  Combining these two statements, we have 
  \[
  	\HOM[_{\inv K}]0{W,V}\subseteq\HOM[_{\inv M}]0{W,V^{\inv N}},
  \]
  for any irreducible $(\inv K)$-representation $W$. Since $V^\inv N$ is $(\mf m, M_0)$-irreducible, the right-hand side has dimension $0$ or $1$, and if it has dimension $1$, then $W\cong V^{\inv N}$ as $(\inv M)$-representations. 

  Assume there is a non-zero $(\inv K)$-invariant vector in $V$. Then there is also a non-zero homogeneous $(\inv K)$-invariant vector $v$ (say). Applying the above to $W=\cplxs v$, we get $\cplxs v\cong V^{\inv N}$ as $(\inv M)$-representations. In particular, $\dim V^{\inv N}=1|0$ or $0|1$, and the parity fixes the parity of $v$. That is, there are no non-zero $(\inv K)$-invariant vectors of the opposite parity. Moreover, $\dim V_j=\dim\HOM[_{\inv K}]0{\cplxs v,V}\sle 1$, where we set $j\defi\Abs0v$, so that $V^{\inv K}=\cplxs v$. This proves the claim. 
\end{PRF}

\subsection{Induced representations}
\label{sec:InducedReps}

In the following, we will consider induced representations. To that end, we first let $(\ger g,G_0)$ be an arbitrary \emph{cs} supergroup pair.

\begin{DFN}
  Let  $(\mf{g}, G_0)$ be a supergroup pair and $W$ a locally convex super vector space. Denote by $\Uenv0{\ger g}$ the universal enveloping algebra of $\ger g$. 

  We define $\Ct[^\infty]0{\ger g,G_0,W}$ to be the set of all linear maps $f:\Uenv0{\ger g}\to\Ct[^\infty]0{G_0,W}$, which are subject to the condition
  \begin{equation}\label{EQ:IndRepInv1}
    \forall x\in\Uenv0{\ger g},y\in\ger g_\ev:f(yx)=\sh L_yf(x).
  \end{equation}
  Here, $\sh L_yf(x)(g)\defi \partial_t|_0f(x)(ge^{ty})$ for $y\in\ger g_{\ev,\reals}$, and this is extended to $\ger g_\ev$ by complex linearity. We call any such $f$ a \emph{$W$-valued superfunction} on $(\mf g,G_0)$.

  We define an action $(\rho_{\mf{g}},\rho_0)$ of $(\mf{g}, G_0)$ on this space by
  \begin{align}
    \big( \rho_0(g) f \big)(u)(p) &\defi f(u)(g^{-1}p),\\
    \big( \rho_{\mf{g}}(x) f \big)(u)(p) &\defi  -(-1)^{|f||x|} f\Parens1{\Ad(p^{-1})(x) u}(p),
  \end{align}
  for all $g,p\in G_0$, $u \in \mf{U}(\mf{g})$, $x \in \mf{g}$.

  The topology on $\Ct[^\infty]0{\ger g,G_0,W}$ is defined as follows. Let $\Norm0\cdot^W_j$, $j\in J$, be a generating set of seminorms on $W$. For any $j\in J$, $u\in\Uenv0{\ger g}$, and any compact set $K\subseteq G_0$, we let $\Norm0\cdot_{j,u,K}$ be the seminorm defined by
  \[
  ||f||_{j,u,K} \defi  \sup_{p \in K} ||f(u)(p)||^W_j.
  \]
  We consider the locally convex topology generated by these seminorms.

\end{DFN}

\begin{PRO}[cinfty-smoothrep]
  The pair $(\rho_\ger g,\rho_0)$ is a smooth representation of $(\ger g,G_0)$ on the locally convex super vector space $\Ct[^\infty]0{\ger g,G_0,W}$.
\end{PRO}

\begin{LMM}[sup-norm]
  Let $(x_a)_{a=1,\dotsc,p}$ be a basis of $\ger g$ and $(x^a)$ the dual basis. Let $j\in J$, $u\in\Uenv0{\ger g}$, $K\subseteq G_0$ be compact, and $L\subseteq\ger g$ be compact. Then
  \[
  \sup_{x\in L}\,\Norm0\cdot_{j,xu,K}\sle\sum_{a=1}^p\sup\,\Abs0{x^a(L)}\,\Norm0\cdot_{j,x_au,K}
  \]
\end{LMM}

\begin{PRF}
  For any $p\in K$ and $x\in L$,
  \[
  f(xu)(p)=\sum_{a=1}^px^a(x)f(x_au)(p).
  \]
  The claim follows immediately from the triangle inequality.
\end{PRF}

\begin{PRF}[\protect{Proof of \ref{pro:cinfty-smoothrep}}]
  We abbreviate $V \defi \Ct[^\infty]0{\ger g,G_0,W}$. It is obvious that $\rho_0$ is an action of $G_0$. Concerning the action of $\ger g$, we introduce the principal anti-automorphism $S$ of $\Uenv0{\ger g}$. This is the linear map $S:\Uenv0{\ger g}\to\Uenv0{\ger g}$ determined by 
  \[
    S(1)=1\ ,\ S(x)=-x\ ,\ S(uv)=(-1)^{\Abs0u\Abs0v}S(v)S(u)
  \]
  \fa $x\in\ger g$ and all homogeneous $u,v\in\Uenv0{\ger g}$. We compute, for $p \in G_0$ $x,y,z \in \mf{U}(\mf{g})$,
  \begin{align}
    \big( \rho_{\mf{g}}(xy) f \big)(z)(p) &= (-1)^{|f|(|x| + |y|) } f\Parens1{ \Ad(p^{-1})\big( S(xy) \big) z}(p) \\
    &= (-1)^{|f|(|x| + |y|) +|x| |y|} f \Parens1{\Ad(p^{-1})\big(S(y)\big) \Ad(p^{-1})\big(S(x)\big) z}(p)
    \\
    &= (-1)^{|f||y|} \big(
    \rho_{\mf{g}}(x) f \big)( \Ad(p^{-1})\big(S(y)\big) z )(p)
    \\
    &=\big( \rho_{\mf{g}}(x) \rho_{\mf{g}}(y) f \big)(z)(p),
  \end{align}
  where we have extended $\rho_{\mf{g}}$ to $\mf{U}(\mf{g})$. It
  follows immediately that $\rho_\mf g$ is a super
  Lie algebra representation.

  Next, we need to check that the $\mf{g}_\ev$-invariance in
  Equation  \eqref{EQ:IndRepInv1} is preserved.  Beginning with $\rho_\ger g$, take $f \in V$, $x \in \mf{U}(\mf{g})$, $y\in \mf{g}_{\ev,\reals}$, $g \in G_0$, and $ v \in \mf{g}$ to compute
  \begin{align}
    (-1)^{|f||v|} &(\rho_{\mf{g}}(v)f)(yx))(g)= f(\Ad(g^{-1})(S(v)) yx)(g)
    \\
    & = - f\Parens1{[\Ad(g^{-1})(S(v)), y] x}(g) + f\Parens1{y\Ad(g^{-1})(S(v))x}(g)
    \\
    &=
    \partial_t \big |_0 \Bracks1{ f\Parens1{\Ad((ge^{ty})^{-1})(S(v)) x}(g) +
      f\Parens1{\Ad(g^{-1})(S(v)) x)(g e^{ty}}}
    \\
    &= (-1)^{|f||v|} \sh L_y((\rho_{\mf{g}}(v)f)(x))(g)
  \end{align}
  where in last step, we have used the product rule. For $\rho_{0}$, the computation is straightforward.

  So far, we have established that $(\rho_\ger g,\rho_0)$ is a well-defined pair of morphisms $\rho_0:G_0\to\Gl(V)$ and $\rho_\ger g:\ger g\to\ger{gl}(V)$ of groups resp.~of of Lie superalgebras. Next, we need to check the conditions of Definition  \ref{DEF:LCSreps}.
  We begin with the continuity of $\rho_0$ and $\rho_\mf g$.

  Let $(g_n, f_n)\in G_0 \times V$ be a net converging to $(g,f)$ and $(x_n,f_n) \in \mf  g \times V$ be a net converging to $(x,f)$. Writing
  \begin{equation}
  	\begin{split}
	    \rho_0(g_n)f_n - \rho_0(g)f &= \rho_0(g) \left(f_n -f\right) + \big(\rho_0(g_n) -\rho_0(g)\big)(f_n),\\
	    \rho_\mf g (x_n) f_n - \rho_\mf g(x)f &=
	    \rho_\mf g (x)( f_n - f) + (\rho_\mf g (x_n) - \rho_\mf g(x))(f_n),  
  	\end{split}
    \label{eq:continuity-2}  	
  \end{equation}
  we can show convergence of each net in two steps.
  Let in the following $j\in J$, $u\in\Uenv0{\ger g}$, and $K\subseteq G_0$ be compact.

  Using \REF{lmm}{sup-norm}, we have
  \[
  \Norm1{\rho_\mf g (x)( f_n - f)}_{j,u,K}
  \leq \sup_{p\in K}\Norm1{ f_n - f}_{j,-\Ad(p^{-1})(x) u,K}
  \to 0
  \]
  since $\Ad(K^{-1})(x)u$ is compact.
  For the second term in Equation \eqref{eq:continuity-2} we assume w.l.o.g.\ that $\rho_\mf g (x)=0$. Then
  we can again use \REF{lmm}{sup-norm} to see that
  \begin{align}
    || \rho_{\mf{g}}(x_n)f_n||_{j,u,K} \leq \sup_{p\in K}
    || f_n||_{j, -\Ad(p^{-1})(x_n) u,K} \to 0.
  \end{align}
  Next we consider $\rho_0(g) (f_n - f)$. By definition
  \[
  \Norm1{\rho_0(g) (f_n - f)}_{j,u,K} = \Norm1{f_n -
    f}_{j,u,g^{-1}K} \to 0.
  \]

  Showing that $(\rho_0(g_n)-\rho_0(g))f_n\to 0$ is a little more tricky. Let $C$ be a compact convex neighbourhood of $0$ in $\ger g_{\ev,\reals}$. For $x\in C$, $g,p\in G_0$, $u\in\Uenv0{\ger g}$, and $f\in V$, we have
  \begin{align}
    f(u)(g^{-1}pe^x)-f(u)(g^{-1}p)=\int_0^1\partial_t f(u)(g^{-1}pe^{tx})\,dt=\int_0^1f(xu)(g^{-1}pe^{tx})\,dt,
  \end{align}
  because $\partial_th(t)=\partial_s|_{s=0}h(t+s)$. Hence,
  \begin{equation}\label{eq:est1}
    \Norm1{f(u)(g^{-1}pe^x)-f(u)(g^{-1}p)}^W_j\sle\Norm0f_{j,xu,g^{-1}p\exp(C)}.
  \end{equation}

  Let $U\subseteq G_0$ be neighbourhoods of $1$ \scth
  $p^{-1}Up\subseteq\exp(C)$ for all $p\in K$. For some $n_0$
  and all $n\sge n_0$, we have $gg_n^{-1}\in U$, so there are
  $x_n(p)\in C$ such that
  $g_n^{-1}p=g^{-1}pe^{x_n(p)}$. Applying
  Equation~\eqref{eq:est1}, we get
  \begin{align}
    \Norm1{(\rho_0(g_n)-\rho_0(g))f_n}_{j,u,K}&=\sup_{p\in K}\Norm1{f_n(u)(g^{-1}pe^{x_n(p)})-f_n(u)(g^{-1}p)}_j^W\\
    &\sle\sup_{x\in C}\Norm0{f_n}_{j,xu,g^{-1}K\exp(C)}
  \end{align}

  Since $C$ and $g^{-1}K\exp(C)$ are compact, and in view of
  Lemma \ref{lmm:sup-norm}, the right-hand side tends to zero
  with $n$. This completes the proof of continuity.

  Next we turn to the relation
  of the two representations.
  %
  For $x \in \mf{g}_{\ev,\R}$, $f \in V$, $u \in \mf{U}(\mf{g})$,
  $p \in G_0$, we have
  \begin{align}
    \partial_t \big |_0 f(u)(e^{-t x}p) &= \partial_t \big |_0
    f(u)(p e^{-t \Ad(p^{-1})x})
    \\
    &=\mathcal{L}_{\Ad(p^{-1})S(x)}f(u)(p)
    \\
    &= f(x u)(p) = \left(\rho_{\mf{g}}(x) f\right) (u)(p),
  \end{align}
  hence $\rho_{\mf{g}}$ extends the derivative of $\rho_0$.
  To verify the intertwining property take $x \in \mf{g}$, $f \in V$, $u \in \mf{U}(\mf{g})$, $p, g
  \in G_0$. Then
  \begin{align}
    \left( \rho_0(g) \rho_{\mf{g}}(x) \rho_0(g^{-1})
      f\right)(u)(p) &=\left( \rho_0(g) \rho_{\mf{g}}(x)
      f\right)(u)(gp)
    \\
    &= (-1)^{|f||x|}\left( \rho_0(g)
      f\right)(\Ad((gp)^{-1})(S(x))u)(gp)
    \\
    &= (-1)^{|f||x|}f(\Ad(p^{-1})\Ad(g)(S(x))u)(p)
    \\
    &= (-1)^{|f||x|}f(\Ad(p^{-1}) (S(\Ad(g)((x)))) u)(p)
    \\
    &= \left( \rho_{\mf{g}}(\Ad(g)x) f\right)(u)(p),
  \end{align}
  where the fact that $S$ commutes with the adjoint action was applied.

  Let $f \in V$. We claim that $f$ is a smooth vector. To that end, let $U\subseteq\ger g_{\ev,\reals}$ be an open neighbourhood of $0$. By the definition of the topology on $V$, it is sufficient to show for any $u\in\Uenv0{\ger g}$ and any $p\in G_0$ the map
  \[
  F:U\to\Ct[^\infty]0{U,W}:x\mapsto\Bracks1{y\mapsto(\vrho_0(e^x)f)(u)(e^yp)}
  \]
  is well-defined and smooth. The map $\phi:U\times U\to G_0:(x,y)\mapsto e^{-x}e^y$ is smooth, and we have $F(x)(y)=f(u)(\phi(x,y)p)$, so this follows from Lemma~\ref{lmm:curry}. All in all, $V$ carries indeed a smooth representation of $(\mf{g},G_0)$.
\end{PRF}

\begin{DFN}[IndRep]
  Let $(\mf{h}, H_0)$ be a sub-supergroup pair of $(\mf{g}, G_0)$ and $\pi=(\pi_\mf{h},\pi_0)$ a smooth representation of $(\mf{h}, H_0)$ on the locally convex super-vector space $W$.  We define $\IND[_{\ger h,H_0}^{\ger g,G_0}]0W$ to be the subspace of $\Ct[^\infty]0{\ger g,G_0,W}$ consisting of all $f$ such that
  \begin{align}
    f(x)(p)&=\pi_0(h)\Parens1{( f(\Ad(h^{-1})x))(ph)}\label{EQ:IndRepInv2}
    \\
    f(xy)(p) &= (-1)^{|y|(|x|+|f|)}
    \pi_{\mf{h}}(S(y)) \big( f(x)(p) \big)
    \label{EQ:IndRepInv3}
  \end{align}
  for all $x\in\Uenv0{\ger g}$, $p\in G_0$, $y\in\mf{h}$, and $h \in H_0$.
\end{DFN}


\begin{LMM}\label{LMM:InducedRepIsRep}
  The subspace $\IND[_{\ger h,H_0}^{\ger g,G_0}]0W$ of $\Ct[^\infty]0{\ger g,G_0,W}$ is invariant under the action of $(\mf{g}, G_0)$. Hence, it defines a smooth representation.
\end{LMM}

\begin{PRF}
  To see that the $H_0$-invariance  
  is  preserved by $\rho_\mf g$ take $h \in H_0$,
  $f \in \IND[_{\ger h,H_0}^{\ger g,G_0}]0W$, $u \in \mf{U}(\mf{g})$,
  $p \in G_0$ and $ x \in \mf{g}$ and
  compute
  \begin{align}
    (\rho_{\mf{g}}(x)f)(u)(p) &= (-1)^{|f||x|} f(\Ad(p^{-1})(S(x))u)(p)
    \\
    &= (-1)^{|f||x|} \pi_0(h) \Big( f(\Ad(h^{-1})(\Ad(p^{-1})(S(x))
    u))(p h) \Big)
    \\
    &= (-1)^{|f||x|} \pi_0(h) \Big( f( \Ad((ph)^{-1})(S(x))
    \Ad(h^{-1})(u))(p h) \Big)
    \\
    &= \pi_0(h) \Big( (\rho_{\mf{g}}(x)f)(\Ad(h^{-1}u)(ph) \Big).
  \end{align}

  It follows immediately from the definitions that
  $\rho_\mf g$ preserves \eqref{EQ:IndRepInv3} and that the action
  $\rho_0$ of $G_0$ preserves \eqref{EQ:IndRepInv2} and
  \eqref{EQ:IndRepInv3}.
  %
  %
\end{PRF}

Now let $\pair g$ be as in the statement of \REF{pro}{mainProof1} and $\pair q$ the minimal parabolic subpair associated with the Iwasawa decomposition.

\begin{PRO}[inducedQGMK]
  For any smooth $(\mf g, G_0)$-representation $(\pi_\mf g, G_0)$ on a locally convex super-vector space $V$, there is an isomorphism of smooth $(\mf k, K_0)$-representations,
  \[
  \Ind^{\mf g, G_0}_{\mf q, Q_0} \left(  V\right)
  \simeq
  \Ind^{\mf k, K_0}_{\mf m, M_0} \left( V \right).
  \]
\end{PRO}

\begin{PRF}
  \allowdisplaybreaks
  Let $\Psi:\Ind^{\ger g,G_0}_{\ger q,Q_0}(V)\to\Ind^{\ger k,K_0}_{\ger m,M_0}(V)$ be the restriction to $\mf{U}(\mf{k}) \subseteq \mf{U}(\mf{g})$ and
  $K_0 \subseteq G_0$. Then $\Psi$ is an injective
  morphism of $(\mf k, K_0)$ representations, because
  $\mf{U}(\mf{g})\simeq\mf{U}(\mf{k})\otimes \mf{U}(\mf{a} \oplus \mf{n})$ as
  vector spaces and
  \[
  f(xy)(kan)=(-1)^{(|f| + |x|)(|y|)} \pi_{\mf{g}}(S(y)) \pi_0(an)^{-1} f(x)(k)
  \]
  for $x \in \mf{U}(\mf{k})$, $y \in \mf{U}(\mf{a}\oplus\mf{n})$,$k \in K_0$, $a
  \in A_0$ and  $n \in N_0$.

  To define the inverse we choose a homogeneous basis $(e_i)$ of $\mf{U}(\mf{a} \oplus \mf{n})$ and denote
  \begin{align}
    \mf{U}(\mf{g}) \xrightarrow{\sim} \mf{U}(\mf{k})\otimes \mf{U}(\mf{a} \oplus
    \mf{n})
    \label{EQ:Ug=UkUan}
    \;:\;
    x \mapsto \sum_i [x]_i \otimes e_i
  \end{align}
  where the sum is finite for all $x$.
  We define
  $\Phi:\Ind^{\mf k, K_0}_{\mf m, M_0} \left( V \right)
  \to\Ind^{\mf g, G_0}_{\mf q, Q_0} \left(  V\right) $
  by
  \[
  \Phi(f)(x)(kan)\defi \sum_i \pm\pi_0(an)^{-1}\pi_{\mf{g}}(S(e_i))  f([\Ad(an)(x)]_{i})(k)
  \]
  where the sign is given by $(-1)^{(|f| +  |[\Ad(an)(x)]_{i}|)|e_i|}$.

  Formally, it is obvious that $\Psi \circ \Phi=\id$, and an easy computation shows $\Phi\circ\Psi=\id$, but to see that $\Phi$ is well-defined we need to check the Conditions \eqref{EQ:IndRepInv1}, \eqref{EQ:IndRepInv2}, and \eqref{EQ:IndRepInv3}. To that end, let $g=kan\in G_0=K_0A_0N_0$ and $p=m' a' n' \in Q_0 = M_0 A_0 N_0$. Then
  \[
  gp=kan m'a'n' = \left( k m' \right) \left( a a' \right) \left( (m'
    a')^{-1} n (m' a') n' \right) \in K_0 A_0 N_0
  \]
  where we have used the fact that $M_0$ centralises $A_0$ and $M_0A_0$ normalises $N_0$.

  As a shorthand, we set, for $x \in \mf{U}(\mf{g})$,
  \[
  z \defi  \Ad\left( \left( a a' \right) \left( (m'
      a')^{-1} n (m' a') n' \right)\right)\left(\Ad(m'a'n')^{-1}(x)\right)
  =\Ad\left( a (m')^{-1} n \right)(x).
  \]
  For $f \in \Ind_{\mf m, M_0}^{\mf k, K_0}(V)$, we compute
  \begin{subequations}
    \begin{align}
      &\pi_0(p) \left(
        \Phi(f)(\Ad(p)^{-1}(x))(gp)
      \right)
      \\
      &=
      \pi_0(p)
      \sum_i
      (-1)^{(|f|+|\left[ z \right]_{i}|)(|e_i|)}
      \pi_0(a (m')^{-1} n m'a'n')^{-1}
      \pi_{\mf{g}}(S( e_i))
      f( [z]_{i} )(km')
      \\
      &=\sum_i
      (-1)^{(|f|+|\left[ z \right]_{i}|)(|e_i|)}
      \pi_0(a (m')^{-1} n)^{-1}
      \pi_{\mf{g}}(S( e_i))
      \pi_0(m')^{-1}
      f( \Ad(m')([z]_{i}) )(k)
      \label{EQ:RestrictInd-1}
      \\
      &=\sum_i
      (-1)^{(|f|+|[\Ad(an)(x)]_{i}|)(|e_i|)}
      \pi_0(a n)^{-1}
      \pi_{\mf{g}}(S( e_i))
      f( [\Ad(an)(x)]_{i}) )(k)
      \\
      &=
      \Phi(f)(x)(g)
    \end{align}
  \end{subequations}
  where we have employed the $M$-equivariance of $f$, and that  $\Ad(M_0)$ preserves the decomposition
  \eqref{EQ:Ug=UkUan} and hence commutes with the projections. We
  have verified \eqref{EQ:IndRepInv2}.

  To check equivariance \eqref{EQ:IndRepInv3} at the algebra level, let $y \in
  \mf{U}(\mf{q})$. Further, denote
  $K_i \defi  [\Ad(an)x]_{i} $
  and $M_j \defi  [\Ad(an)y]_{j} $. Then
  \[
  K_ie_iM_je_j=K_i[e_i,M_j]e_j + (-1)^{|e_i||M_j|}K_iM_je_ie_j
  \]
  with $[e_i,M_j], e_ie_j \in \mf{U}(\mf{a} \oplus \mf{n})$ and $K_i, K_iM_j \in \mf{U}(\mf{k})$. We compute
  \begin{align}
    &\pi_0(an)\Phi(f)(xy)(g)
    \\
    &=\sum_{i,j}\begin{aligned}[t]
      &(-1)^{(|f|+|K_i|)(|e_i|+|M_j|+|e_j|)}\pi_{\mf{g}}\left( S([e_i,M_j]e_j) \right)f(K_i)(k)\\
      &+(-1)^{(|f|+|K_i|+|M_j|)(|e_i|+|e_j|) +|e_i||M_j|}
      \pi_{\mf{g}}\left( S(e_ie_j) \right)f(K_iM_j)(k)
    \end{aligned}\\
    &=\sum_{i,j}
    \begin{aligned}[t]
      &(-1)^{(|f|+|K_i|)(|e_i|+|M_j|+|e_j|)}\pi_{\mf{g}}\left( S(e_iM_je_j) \right)f(K_i)(k)\\
      &-(-1)^{(|f|+|K_i|)(|e_i|+|M_j|+|e_j|) +|e_i||M_j|}
      \pi_{\mf{g}}\left( S(M_je_ie_j) \right)f(K_i)(k)\\
      &+(-1)^{(|f|+|K_i|)(|e_i|+|e_j|+|M_j|)+|M_j||e_j|}
      \pi_{\mf{g}}\left( S(e_ie_j) \right)\pi_{\mf{g}}(S(M_j))f(K_i)(k)
    \end{aligned}\\
    &=
    \sum_i(-1)^{(|f|+|K_i|)(|e_i|+|y|)+|e_i||y|}\pi_0(an)\pi_{\mf{g}}(S(y))
    \pi_0(an)^{-1}
    \pi_{\mf{g}}\left( S(e_i) \right)
    f(K_i)(k)
    \\
    &=
    (-1)^{(|f|+|x|)|y|}\pi_0(an)\pi_{\mf{g}}(S(y))
    \Phi(f)(x)(g),
  \end{align}
  hence $\Phi(f)$ is $\mf q$-equivariant.

  To check Condition \eqref{EQ:IndRepInv1}, let $x \in\mf{g}_{\ev,\reals}$ and denote by $k_t$ and $s_t$, respectively, the $K_0$ and $A N_0$ part of $e^{ t \Ad(an) x}$,  for small $t \in \R$. This defines curves in $K_0$ and $AN_0$ whose derivatives at $t=0$ we denote by $\dot k_0$ and $\dot s_0$.  Furthermore, we let $M_{jt}\defi  [\Ad(s_tan)y]_j$. These are curves in a finite-dimensional subspace of $\mf{U}(\mf{k})$. We denote their derivatives at $t=0$ by $\dot{M}_{j0}$.
  Then
  \begin{align}
    &\partial_t \big |_0 \Phi(f)(y)(ge^{t x})
    =\partial_t \big |_0 \Phi(f)(y)(k e^{ t \Ad(an) x} an)
    \\
    &=
    \sum_j\partial_t \big |_0
    (-1)^{(|f|+|y|+|e_j])|e_j|}
    \pi_0(an)^{-1}\pi_0(s_t)^{-1}\pi_{\g}(S(e_j)) f(M_{jt})(kk_t)
    \\
    &=\sum_j
    \begin{aligned}[t]
      &(-1)^{(|f|+|y|+|e_j|)|e_j|} \pi_0(an)^{-1} \big(\pi_{\g}(S(e_j\dot s_0 )) f(M_{j0})(k)\\
      &+\pi_{\g}(S(e_j)) f( \dot{M}_{j0})(k)+\pi_{\g}(S(e_j)) f( \dot{k}_0 M_{j0})(k)\big)
    \end{aligned}
    \\
    &=
    \sum_j\Phi(f) \big( \Ad(an)^{-1} \big(
    M_{j0} e_j \dot s_0
    + \dot{M_{j0}} e_j
    +\dot{k}_0 M_{j0} e_j \big)\big)(g)
    \\
    &=
    \Phi(f) \big( y \Ad(an)^{-1} (\dot s_0 )
    +[\Ad(an)^{-1} \dot s_0, y]
    + \Ad(an)^{-1}(\dot{k}_0) y \big)(g)
    \\
    &=
    \Phi(f) \big( \Ad(an)^{-1} ( \dot{k}_0+\dot s_0) y \big)(g)=\Phi(f) \big( x y \big)(g)
  \end{align}
  The required smoothness of $\Phi(f)$ is immediate from its definition. So, all in all, $\Phi$ is well-defined, and it is certainly an even continuous linear map.

  Since $\Phi$ is inverse to $\Psi$, it is $\pair k$-equivariant, so $\Psi$ is an isomorphism.
\end{PRF}
%


\begin{PRO}[embedding]
  Let $(\pi_0, \pi_{\g})$ an irreducible finite-dimensional representation of $(\mf g, G_0)$ on $V$. There is an injective morphism of smooth $\pair g$-representations
  \[
  \alpha: V \into \Ind_{\mf q, Q_0}^{\mf g, G_0}(V^\inv N).
  \]
\end{PRO}

The proof is divided into the following \REF{lmm}{indDoubleDual} and \REF{lmm}{embedIntoDD}.

\begin{LMM}[indDoubleDual]
  There is an isomorphism of smooth $\pair Q$-representations
  \[
  V^\inv N\simeq\Parens1{\Parens0{  V^* }^{\inv{\bar{N}}}}^*.
  \]
\end{LMM}

\begin{PRF}
  Let
  \[
  \Phi: V^\inv N\to\Parens1{(   V^* )^{\inv{\bar{N}}}}^*, \Phi(v)(\mu)\defi \mu(v).
  \]
  We will show that $\Phi$ is an isomorphism of $(\mf q,Q_0)$-representations. By definition of the dual representation (\REF{dfn}{dualRep}), $\Phi$ is $(\mf q,
  Q_0)$-equivariant.

  Further we observe that $\Phi$ is injective. Indeed, let $\Phi(v)=0$ for some $v\in V^\inv N$. Let $\mu\in(V^*)^{\inv{\bar N}}$. Then $\mu(v)=0$ and for any $x\in\ger n$,
  \[
  (\pi_{\g}^*(x)\mu)(v)=-\mu(\pi_{\g}(x)(v))=0,
  \]every  since $v$ is $\ger n$-invariant. The lowest weight vector of $V^*$ is contained in $(V^*)^{\inv{\bar{N}}}$, so $\pi_{\g}^*(\Uenv0{\ger n})((V^*)^{\inv{\bar{N}}}) = V^*$. It follows that $\mu(v)=0$ for every $\mu\in V^*$, so $v=0$.

  Hence, $\dim(V^\inv N)^* \leq \dim(V^*)^{\inv{\bar{N}}}$. Interchanging the
  roles of $N$ and $\bar{N}$, the dimensions are equal, and $\Phi$ is an
  isomorphism.
\end{PRF}

\begin{LMM}[embedIntoDD]
  There is an injective morphism of smooth $(\mf g, G_0)$-representations
  \begin{gather}
    \alpha : V\into \Ind^{\mf g, G_0}_{\mf q, Q_0}\Parens1{\Parens0{  V^* }^{\inv{\bar{N}}}}^*,\\
    \alpha(v)(u)(p)(\mu) \defi  (-1)^{|u||v| + |\mu|(|u|+|v|)} \mu \Parens1{
      \pi_0(p^{-1}) \pi_{\mf{g}}(\Ad(p)S(u))v }.
  \end{gather}
\end{LMM}

\begin{PRF}

  Slightly abusing notation, we will denote by $(\pi_\ger g^{**},\pi_0^{**})$ the $\pair q$-action on $((V^*)^{\inv{\bar n}})^*$. The map
  $\alpha$ is well-defined. Indeed, for $p \in G_0$, $v \in V$, $u \in \mf{U}(\mf{g})$, $\mu \in (V^*)^{\inv{\bar{N}}}$, and $q \in Q_0$,
  \begin{multline}
    \Bracks1{\pi_0^{**}(q)\Parens1{\alpha(v)(\Ad(q^{-1})u)(pq)}}(\mu)
    \\
    =
    \pm\mu \Parens1{ \pi_0(q) \pi_0((pq)^{-1}) \pi_{\mf{g}}(\Ad(pq)S(\Ad(q^{-1})u)) v
    }
    =\alpha(v)(u)(p)(\mu).
  \end{multline}
  Similarly, for $y \in \mf{q}$,
  \begin{align}
    \alpha(v)(uy)(p)(\mu)
    &=
    (-1)^{(|u|+|y|)|v| + |\mu|(|u|+|y|+|v|)}\mu \left(
      \pi_{\mf{g}} (S(uy))\pi_0(p^{-1}) v
    \right)
    \\
    &=
    (-1)^{(|u|+|v|)|y| }\Parens1{\pi_{\mf{g}}^{**}(S(y)) \alpha(v)(u)(p)} (\mu),
  \end{align}
  and for $x \in \mf{g}$,
  \begin{align}
    \alpha(v)(x u)(g)(\mu)
    &=
    (-1)^{|u||v| + |\mu|(|u|+|v|)}\mu \left(
      \pi_{\mf{g}} (S(x u))\pi_0(g^{-1}) v
    \right)
    \\
    &=
    (-1)^{|u||v| + |\mu|(|u|+|v|)}\mu \left(
      \pi_{\mf{g}} (S(u)) \pi_{\mf{g}}(-x)\pi_0(g^{-1}) v
    \right)
    \\
    &=
    (-1)^{|u||v| + |\mu|(|u|+|v|)}\mu \left(
      \pi_{\mf{g}} (S(u))  \partial_t \big|_0  \pi_0(e^{-t x}g^{-1}) v
    \right)
    \\
    &=
    \partial_t \big|_0 \alpha(v)(u)(g e^{t x})(\mu).
  \end{align}

  Smoothness of $\alpha(v)(u)$ follows from the smoothness of $\pi_0$ and finite-dimensiona-lity of $V$. Hence, $\alpha(v)$
  indeed lies in the induced representation.

  Further, $\alpha$ is a morphism of smooth $\pair G$-representations. Then,
  for $g \in G_0$,
  \begin{align}
    \alpha(\pi_0(g)v)(u)(p)(\mu) &= (-1)^{|u||v| + |\mu|(|u|+|v|)} \mu \left(
      \pi_{\mf{g}}(S(u))  \pi_0(p^{-1})  \pi_0(g) v \right)
    \\
    &= (-1)^{|u||v| + |\mu|(|u|+|v|)} \mu \left(
      \pi_{\mf{g}}(S(u))  \pi_0((g^{-1}p)^{-1}) v \right)
    \\
    &=(\rho_0(g)\alpha(v))(u)(p)(\mu).
  \end{align}
  Similarly, for $x \in \mf{g}$,
  \begin{align}
    &\alpha(\pi_{\g}(x)v)(u)(p)(\mu)
    \\
    &= (-1)^{|u|(|v|+|x|) + |\mu|(|u|+|x|+|v|)}
    \mu \left(  \pi_0(p^{-1}) \pi_{\mf{g}}(\Ad(p)S(u))  \pi_{\g}(x)   v \right)
    \\
    &= (-1)^{|u||v| + |\mu|(|u| +|x|+|v|)}
    \mu \left(
      \pi_0(p^{-1}) \pi_{\mf{g}}(\Ad(p)S\left( \Ad(p^{-1})(S(x)) u\right ))  v \right)
    \\
    &= (-1)^{|x||v| }\alpha(v)(\Ad(p^{-1})(S(x))u)(p)(\mu)= (\rho_{\g}(x)\alpha(v))(u)(p)(\mu).
  \end{align}

  Finally, $\alpha \neq 0$ because
  \[
  \alpha(v_{\lambda})(1)(\1)(\mu_{-\lambda}) =
  \mu_{-\lambda}(v_{\lambda}) \neq 0
  \]
  for the highest weight vector
  $v_{\lambda}$ of $V$ and the lowest weight vector $\mu_{-\lambda}$ of $V^*$. Since $V$ is irreducible, $\alpha$ is injective.
\end{PRF}

\subsection{Frobenius reciprocity}
\label{sec:frobenius}

In proving \REF{pro}{mainProof1}, we have applied a super version of Frobenius reciprocity. To state it in full generality, let $(\ger g,G_0)$ be an arbitrary \emph{cs} supergroup pair.

\begin{PRO}[Frobenius]
  Let $(\mf{g},G_0)$ be a supergroup pair with subgroup pair $(\mf{h},H_0)$. For any smooth representation $V$ of $(\ger g,G_0)$ and any smooth representation $W$ of $(\ger h,H_0)$, we have a natural isomorphism
  \[
  \HOM[_{\mf{g},G_0}]1{V, \Ind^{\mf{g},G_0}_{\mf{h},H_0}(W)} \simeq \HOM[_{\mf{h},H_0}]0{V, W}.
  \]
  On the right-hand side, $V$ is considered as a representation of $(\mf{h},H_0)$.
\end{PRO}

\begin{PRF}
  Consider the linear map
  \begin{align}
    \Phi : \Hom_{\mf{g},G_0}(V, \Ind^{\mf{g},G_0}_{\mf{h},H_0}(W)) &\to \Hom_{\mf{h},H_0}(V, W), \Phi(T)(v)\defi T(v)(1)(\1),
  \end{align}
  where $1 \in \mf{U}(\mf{g})$ is the multiplicative unit and $\1 \in G_0$ is the neutral element.

  To see that $\Phi$ is well-defined take $T \in \Hom_{\mf{g},G_0}(V, \Ind^{\mf{g},G_0}_{\mf{h},H_0}(W))$, $v\in V$, $ h \in H_0$, and $y \in \mf{h}$, and compute
  \begin{subequations}
    \begin{align}
      \Phi(T)(\pi^V_0(h)v) &= T(\pi^V_0(h)v)(1)(\1)= (\rho_0(h)T(v))(1)(\1) \label{EQ:prfFrobenius-1}\\
      &= T(v)(1)(h^{-1})= \pi^W_0(h) \big(T(v)(1)(\1)\big)\label{EQ:prfFrobenius-2}\\
      &=\pi^W_0(h)\Parens1{\Phi(T)(v)}
    \end{align}
  \end{subequations}
  where in Equation~\eqref{EQ:prfFrobenius-1}, we used the $G_0$-equivariance of $T$, and in Equation
  \eqref{EQ:prfFrobenius-2}, we used the fact that $T(v)$ is a vector in the induced representation.  Similarly,
  \begin{align}
    \Phi(T)(\pi^V_{\mf{g}}(y)v) &= (\rho_{\mf{g}}(y)T(v))(1)(\1)= (-1)^{|T(v)| |y|} \;T(v)(S(y))(\1)\\
    &= (-1)^{|T(v)| |y|} (-1)^{(|T(v)|+|1|) |y|} \;\pi^W_{\mf{h}}(y)\Parens1{T(v)(1)(\1)}
    \\
    &=\pi^W_{\mf{h}}(y)\Parens1{\Phi(T)(v)}
  \end{align}
  The continuity of $T$ implies continuity of $\Phi(T)$.  Hence, $\Phi$ is well defined.

  To see that $\Phi$ is invertible, let $u \in \mf{U}(\mf{g})$ and $g \in G_0$ and notice
  \begin{align}
    T(v)(u)(g) &= (-1)^{|v||u|} \left( \rho_{\mf{g}}(S(u)) \circ
      \rho_0(g^{-1}) T(v) \right)(1)(\1)
    \\
    &= (-1)^{|v||u|} T\left(\pi_{\mf{g}}^V(S(u)) \pi_0^V(g^{-1}) v
    \right)(1)(\1)
    \\
    &= (-1)^{|v||u|} \Phi(T)\left( \pi_{\mf{g}}^V(S(u))
      \pi_0^V(g^{-1}) v \right)
  \end{align}
  which means that $\Psi\circ \Phi (T) = T$ for
  \begin{align}
    &\Psi: \Hom_{\mf{h},H_0}(V, W) \to \HOM[_{\mf{g},G_0}]1{V,
      \Ind^{\mf{g},G_0}_{\mf{h},H_0}(W)},\\
    &\Psi(R)(v)(u)(g)\defi (-1)^{|v||u|} R\Parens1{ \pi_{\mf{g}}^V(S(u)) \pi_0^V(g^{-1}) v}
  \end{align}

  This definition of $\Psi$ also yields $\Phi\circ\Psi(R)=R$, formally. However, it remains to be shown that $\Psi$ is well-defined in the first place. To that end, take $p\in G_0$ and $u \in \mf{U}(\mf{g})$ and compute
  \begin{align}
    \Psi(R)(\pi_0^V(g)v )(u)(p) &= (-1)^{|v| |u|} R\left(
      \pi_{\mf{g}}^V(S(u)) \pi_0^V(p^{-1}) \pi_0^V(g) v \right)
    \\
    &= (-1)^{|v| |u|} R\left( \pi_{\mf{g}}^V(S(u))
      \pi_0^V((g^{-1}p)^{-1}) v \right)
    \\
    &= \Psi(R)(v )(u)(g^{-1}p) =
    \left(\rho_0(g)\Psi(R)(v)\right)(u)(p).
  \end{align}
  Similarly, with $x \in \mf{g}$,
  \begin{align}
    \Psi(R)(\pi_{\mf{g}}^V(x)v )(u)(p) &= (-1)^{(|v|+|x|) |u|}
    R\left( \pi_{\mf{g}}^V(S(u)) \pi_0^V(p^{-1}) \pi_{\mf{g}}^V(x) v
    \right)
    \\
    &= (-1)^{(|v|+|x|) |u|} R\left( \pi_{\mf{g}}^V(S(u)
      \Ad(p^{-1})(x) ) \pi_0^V(p^{-1}) v \right)
    \\
    &= (-1)^{|v||u|} R\left( \pi_{\mf{g}}^V( S
      (S(\Ad(p^{-1})(x)) u) ) \pi_0^V(p^{-1}) v \right)
    \\
    &= (-1)^{\Abs0v\Abs0x}\Psi(R)(v )(\Ad(p^{-1})(S(x)) u)(p)
    \\
    &= \left( \rho_{\mf{g}}(x)\Psi(R)(v ) \right) (u)(p).
  \end{align}
  Hence, $\Psi(R)$ is indeed $(\ger g,G_0)$-equivariant.

  Next, we need to show that $\Psi(R)(v) \in
  \Ind_{\mf{h},H_0}^{\mf{g},G_0}(W)$.
  Firstly, because $R$ and $\smash{\pi^V_\ger g(S(u))}$ are continuous and linear (hence smooth), and because $v$ is a smooth vector, we find that $\Psi(R)(v)(u)\in\Ct[^\infty]0{G_0,W}$. Moreover, this quantity evidently depends linearly on $u$.

  To check Condition~\eqref{EQ:IndRepInv1} let $y \in \g_\ev$, $u \in \mf{U}(\g)$,
  and $p \in G_0$, and compute
  \begin{align}
    \Psi(R)(v)(y u)(p) &= (-1)^{|v| |u|} R(\pi^V_{\g}(S(y
    u))\pi^V_0(p^{-1})v)
    \\
    &=
    (-1)^{|v| |u|} R(\pi^V_{\g} (S(u)) \pi^V_{\g}(-y)\pi^V_0(p^{-1})v)
    \\
    &=
    \partial_t \big |_0 (-1)^{|v| |u|} R(\pi^V_{\g} (S(u)) \pi^V_0(e^{-t
      y}p^{-1})v)
    \\
    &=
    \partial_t \big |_0 \Psi(R)(v)( u)(p e^{ty}).
  \end{align}
  Thus, we have $\Psi(T)(v)\in\Ct[^\infty]0{\ger g,G_0,W}$.

  For Condition~\eqref{EQ:IndRepInv2}, let $h \in H_0$. Then
  \begin{align}
    \Psi(R)(v)(\Ad(h)(u))(ph^{-1}) &=
    (-1)^{|v| |u|}
    R \big(\pi^V_{\g}(S( \Ad(h)(u))\pi^V_0(hp^{-1})v\big)
    \\
    &= (-1)^{|v| |u|} \pi_0^W(h) R(\pi^V_{\g}(S(u))\pi^V_0(p^{-1})v)
    \\
    &= \pi_0^W(h) \Psi(R)(v)(u)(p)
  \end{align}
  where we used that $R$ is $(\mf{h},H_0)$-equivariant. Finally, to check Condition~\eqref{EQ:IndRepInv3}, let $x \in \mf{h}$. Then
  \begin{align}
    \Psi(R)(v)(ux)(g) &=
    (-1)^{|v| |u|} R(\pi^V_{\g}(S(ux))\pi^V_0(g^{-1})v)
    \\
    &=
    (-1)^{|v| |u| +|u||x|} R( \pi^V_{\g}(S(x)) \pi^V_{\g}(S(u)) \pi^V_0(g^{-1})v)
    \\
    &=
    (-1)^{|u||x|} \pi^W_{\ger h}(S(x))
    \Psi(R)(v)(ux)(g),
  \end{align}
  where we have used the equivariance of $R$ again.

  So $\Psi(R)(v) \in \Ind_{\mf{h},H_0}^{\mf{g},G_0}(W)$. The definition of the topology on $\Ind^{\mf{g},G_0}_{\mf{h},H_0}(W)$ and the continuity of $R$ and $\pi^V$ imply the continuity of $\Psi(R)$. Thus, $\Psi$ is indeed well-defined and inverse to $\Phi$.
\end{PRF}

\section{Representations with $M$-invariant highest weight are
  spherical}
\label{sec:mainProof2}

In this section, we prove necessity in \REF{thm}{main}, stated below as \REF{pro}{mainProof2}. To that end, we will explicitly construct a spherical vector from the highest weight vector in much the same way as in the even case. 
However, deciding whether this vector is zero or not is a more delicate
matter in our super setting. A sufficient criterion for this is furnished by the non-vanishing of the Harish-Chandra $c$-function for the symmetric superspace $G/K$, which we compute explicitly. 


Let $(\mf{g}, G_0)$, $\pair k$, $\pair m$, $\pair a$, and $\pair n$ be
as in the statement of \REF{thm}{main}. We let $\pair q$ denote the
minimal parabolic subpair, defined by $\mf{q}\defi \mf{m}\oplus \mf{a}
\oplus \mf{n}$ and $Q_0\defi M_0A_0N_0$. We also consider the associated Lie supergroups, denoted by $G$, $K$, $M$, $A$, $N$, and $Q$, respectively. 

It what follows, recall the formalism of generalised points and basic facts on Berezin integration, as summarised in Appendix \ref{app:point}.

\begin{PRO}[mainProof2]
  If $V$ is a finite-dimensional irreducible smooth
  $G$-representation whose highest weight $\lambda$ is high enough, then $V$ is spherical.
\end{PRO}

In the \emph{proof}, we take note of the following lemma, which is of separate interest.

\begin{LMM}[knbar-shift]
For $f,h\in\Gamma(\sh O_K)$ and $g\in_SG$, we have 
\[
  \int_Kf(k(g^{-1}k))h(k)e^{\lambda(H(g^{-1}k))}\,\Abs0{Dk}=\int_Kf(k)h(k(gk))e^{-(\lambda+2\vrho)(H(gk))}\,\Abs0{Dk},
\]
where $\vrho\defi\frac12\str_{\ger n}\ad|_{\ger a}=\frac12\sum_{\beta\in\Sigma^+}m_\beta\beta$ and $m_\beta\defi\dim\ger g^\beta_\ev-\dim\ger g^\beta_\odd$. 

\end{LMM}

\begin{PRF}
  Let $\chi\in\Gamma(\sh O_{AN})$ \scth $\int_{AN}\chi(an)e^{2\vrho(\log a)}\,da\,\Abs0{Dn}=1$. We compute, using the Iwasawa decomposition, that 
  \begin{align*}
    \int_K&f(k(g^{-1}k))h(k)e^{\lambda(H(g^{-1}k))}\,\Abs0{Dk}\\
    &=\int_{KAN}f(k(g^{-1}k))h(k)\chi(an)e^{\lambda(H(g^{-1}k))+2\vrho(\log a)}\Abs0{Dk}\,da\,\Abs0{Dn}\\
    &=\int_Gf(k(g^{-1}g'))h(k(g'))\chi(k(g')^{-1}g')e^{\lambda(H(g^{-1}g')-H(g'))}\,\Abs0{Dg'}\\
    &=\int_Gf(k(g'))h(k(gg'))\chi(k(gg')^{-1}g')e^{\lambda(H(g')-H(gg'))}\,\Abs0{Dg'}\\
    \intertext{by \cite{a-hchom}*{Proposition 2.2} and the invariance of $\Abs0{Dg'}$. Furthermore, by the invariance of $da$ and $\Abs0{Dn}$, we find that this equals}
    &=\int_{KAN}f(k)h(k(gk))\chi(e^{H(gk)}n(gk)an)e^{-\lambda(H(gk))+2\vrho(\log a)}\,\Abs0{Dk}\,da\,\Abs0{Dn}\\
    &=\int_{KAN}f(k)h(k(gk))e^{-(\lambda+2\vrho)(H(gk))}\chi(an)e^{2\vrho(\log a)}\,\Abs0{Dk}\,da\,\Abs0{Dn}\\
    &=\int_Kf(k)h(k(gk))e^{-(\lambda+2\vrho)(H(gk))}\,\Abs0{Dk},
  \end{align*}
  where in the last step, the defining property of $\chi$ was applied again.
\end{PRF}

\begin{PRF}[\protect{Proof of \REF{pro}{mainProof2}}]
  Let $\mu_{-\lambda}\neq0$ be a lowest weight vector of $V^*$. Define a linear map
  \[
  \beta:\IND[_{\pair q}^{\pair g}]0{V^{\ger n,N_0}}=(\Gamma(\sh O_G)\otimes V^N)^Q\to V
  \]
  by
  \[
  \beta(f)\defi\int_{K/M}\Dual0{\mu_{-\lambda}}{f(k)}\pi(k)v_\lambda\,\Abs0{D\dot k}.
  \]
  This is well-defined, since for $k\in_SK$, $m\in_SM$, we have
  \[
  \Dual0{v_{-\lambda}^*}{f(km)}\pi(km)v_\lambda=\Dual0{\pi^*(m)v_{-\lambda}^*}{f(k)}\pi(k)v_\lambda=\Dual0{v_{-\lambda}^*}{f(k)}\pi(k)v_\lambda.
  \]
  We observe that $\dim K/M=\dim N=*|2q$ for some $q$, so that $\beta$ is even.

  Next, $\beta$ is $G$-equivariant. Indeed, for $f\in\IND[_{\ger q,Q_0}^{\ger g,G_0}]0{V^N}$, $k\in_SK$, $g\in_SG$,
  \[
  f(g^{-1}k)=f(k(g^{-1}k)e^{H(g^{-1}k)})=e^{-\lambda(H(g^{-1}k))}f(k(g^{-1}k)),
  \]
  and so
  \[
  \beta(gf)=\int_{K/M}e^{-\lambda(H(g^{-1}k))}\Dual0{\mu_{-\lambda}}{f(k(g^{-1}k))}\pi(k)v_\lambda=\Dual0{\mu_{-\lambda}}{\pi(k)f(g^{-1}k)}
  \]
  since we may apply the integral identity in \REF{lmm}{knbar-shift} after rewriting the integrand as an integral over $M$, \cf Ref.~\cite{ah-berezin}*{Corollary 5.12}.

  By assumption, $f\defi1\otimes v_\lambda$ is contained in $\IND[_{\ger q,Q_0}^{\ger g,G_0}]0{V^N}$ and certainly defines a $K$-invariant vector there. Hence,
  \[
  v_K\defi\beta(f)=\Dual0{\mu_{-\lambda}}{v_\lambda}\int_{K/M}\pi(k)v_\lambda\,\Abs0{D\dot k}
  \]
  is a $K$-invariant vector in $V$. We need to see that it is non-zero under the assumption on $\lambda$. It is sufficient to see that
  \[
  \int_{K/M}\Dual0{\mu_{-\lambda}}{\pi(k)v_\lambda}\,\Abs0{D\dot k}\neq0.
  \]
  This is proved in \REF{pro}{km-cfn}.
\end{PRF}

\subsection{Integrals of matrix coefficients}\label{ss:mat}

In this subsection, we reduce the non-vanishing of the integral over $K/M$ considered in the proof of \REF{pro}{mainProof2} to that of an integral over $\bar N$, namely of the Harish-Chandra $c$-function. 

\begin{PRO}[km-cfn]
  We have 
  \begin{equation}\label{eq:kinv-nonzero}
    \int_{K/M}\Dual0{\mu_{-\lambda}}{\pi(k)v_\lambda}\,\Abs0{D\dot k}\neq0,
  \end{equation}
  if and only if $\Dual0{\lambda}{\beta}\neq0$ \fa isotropic positive restricted roots $\beta\in\Sigma^+$, and 
  \[
    \lambda_\alpha+m_\alpha+2m_{2\alpha},\lambda_\alpha+m_\alpha+m_{2\alpha}+1\notin-2\N
  \]
  for all odd anisotropic indivisible restricted roots $\alpha\in\Sigma^+$. In particular, if $\lambda$ is high enough, then the left-hand side of Equation \eqref{eq:kinv-nonzero} is non-vanishing.
\end{PRO}

The \emph{proof} proceeds in several steps, the last of which is \REF{crl}{zeros-cfn}. Firstly, the multiplication map $\bar N\times Q\to G$ is an open embedding, whose underlying open set is dense in $G_0$. We get a corresponding open embedding
\[
\bar N\to G/Q=K/M.
\]
It is given as the composite of the Iwasawa $K$ projection $k:\bar N\to K$ with the canonical morphism $K\to K/M$. So, abusing notation, we write $k$ for this morphism.

\begin{PRO}[knbar-pullback]
  The pullback of the Berezinian density on $K/M$ is given by 
  \[
  k^\sharp\Abs0{D\dot k}=e^{-2\vrho(H(\bar n))}\Abs0{D\bar n},
  \]
  for a suitable normalisation of invariant Berezinian densities.
\end{PRO}

\begin{PRF}
  This is done much as in the ungraded case. Indeed, there exists $\chi\in\Gamma(\sh O_{\bar N})$ \scth $\chi\,\Abs0{D\bar n}=k^\sharp\Abs0{D\dot k}$. Let $f\in\Gamma_c(\sh O_{K/M})$ and $x\in_S\bar N$. Since $k$ is left $AN$-invariant, we have $k(x^{-1}\bar n)=k(x^{-1}k(\bar n))$ for any $\bar n\in_S\bar N$. By invariance of $\Abs0{D\bar n}$ and because $\chi\,\Abs0{D\bar n}=k^\sharp\Abs0{D\dot k}$, we have
  \begin{align*}
    \int_{\bar N} f(k(\bar n))\chi(x\bar n)\,\Abs0{D\bar n}&=\int_{\bar N}f(k(x^{-1}\bar n))\chi(\bar n))\,\Abs0{D\bar n}
    =\int_{K/M}f(k(x^{-1}k))\,\Abs0{D\dot k}.
    \intertext{By \REF{lmm}{knbar-shift}, this equals}
    &=\int_{K/M}f(k)e^{-2\vrho(H(xk))}\,\Abs0{D\dot k}\\
    &=\int_{\bar N}f(k(\bar n))\chi(\bar n)e^{-2\vrho(H(xk(\bar n)))}\,\Abs0{D\bar n}
  \end{align*}
  This gives, for a suitable normalisation, that $\chi(\bar n)=e^{-2\vrho(H(\bar n))}$ for any $\bar n\in_S\bar N$, proving the claim.
\end{PRF}

Now, we apply this coordinate change to the integral in Equation \eqref{eq:kinv-nonzero}. Let $\bar n\in_S\bar N$. Then
\[
\Dual0{\mu_{-\lambda}}{\pi(k(\bar n))v_\lambda}=e^{-\lambda(H(\bar n))}\Dual0{\mu_{-\lambda}}{\pi(\bar n)v_\lambda}
=e^{-\lambda(H(\bar n))}\Dual0{\mu_{-\lambda}}{v_\lambda},
\]
since $\mu_{-\lambda}$ is $\bar N$-invariant. For $f(k)=\Dual0{\mu_{-\lambda}}{\pi(k)v_\lambda}$, \REF{pro}{knbar-pullback} gives
\[
k^\sharp\Parens1{f\,\Abs0{D\dot k}}=\Dual0{\mu_{-\lambda}}{v_\lambda}\,e^{-(\lambda+2\vrho)(H(\bar n))}\,\Abs0{D\bar n}
\]
We are lead to study the integral
\[
c(\lambda+\vrho)\defi\int_{\bar N}e^{-(\lambda+2\vrho)(H(\bar n))}\,\Abs0{D\bar n},
\]
taken with respect to a Weyl retraction of $\bar N$, \vq \REF{dfn}{weyl-ret}. 

We will argue below that this integral converges absolutely, in the sense of Appendix \ref{app:point}, and the outcome $c(\lambda+\vrho)$ is non-zero, provided $\lambda$ is high enough (\REF{crl}{zeros-cfn}). Moreover, we will show in \REF{pro}{rankred} that 
\[
  c(\lambda+\vrho)=\int_{K/M}\frac{\Dual0{\mu_{-\lambda}}{\pi(k)v_\lambda}}{\Dual0{\mu_{-\lambda}}{v_\lambda}}\,\Abs0{D\dot k}.
\]
Thus, the statement of \REF{pro}{km-cfn} follows, thereby completing the proof of \REF{pro}{mainProof2}, and hence, of \REF{thm}{main}.

 \subsection{The $c$-function}

 In this subsection, we prove the convergence of the integral
 $c(\lambda)$ and determine its value.
The general case will be proved
by `rank reduction', much as for the classical Gindikin--Karpelevic
formula. The main difference is that it is not sufficient to consider the action of the Weyl group in this procedure. Rather, it has to be extended to cover also the case of `odd reflections'. The outcome is the following.

\begin{THM}[full-cfn]
  Let $\Re\Dual0\lambda\alpha>0$ for all $\alpha\in\Sigma^+$, $\Dual0\alpha\alpha\neq0$. Then the integral $c(\lambda)$ converges for any Weyl retraction, \vq \REF{dfn}{weyl-ret}, and equals
  \[
  c(\lambda)=
  c_0\prod_{\Dual0\alpha\alpha\neq0}2^{-\lambda_\alpha}\frac{\Gamma(\lambda_\alpha)}{\Gamma\Parens1{\tfrac12(\tfrac{m_\alpha}2+1+\lambda_\alpha)}\Gamma\Parens1{\tfrac12(\tfrac{m_\alpha}2+m_{2\alpha}+\lambda_\alpha)}}\prod_{\Dual0\alpha\alpha=0}\Dual0\lambda\alpha^{-\frac{m_\alpha}2}
  \]
  for some non-zero constant $c_0$, independent of $\lambda$. Here, $\lambda_\alpha\defi\frac{\Dual0\lambda\alpha}{\Dual0\alpha\alpha}$ for $\Dual0\alpha\alpha\neq0$, and the product extends over all indivisible positive restricted roots.
\end{THM}

\begin{PRF}
This follows from the rank reduction formula given in \REF{pro}{rankred} and
the explicit formulae for the anisotropic rank one case, \REF{pro}{cfn-rk1}, and
the isotropic case, \REF{pro}{cfn-isotropic}.
\end{PRF}

From \REF{thm}{full-cfn}, we read off the zeros of the $c$-function. 

\begin{CRL}[zeros-cfn]
We have 
$c(\lambda+\rho)=0$ if and only if one of the following
conditions is fulfilled for some indivisible restricted root $\alpha$:
\begin{enumerate}[wide]
\item $\alpha$ is isotropic and $\Dual0\lambda\alpha=0$.
\item $\alpha$ is anisotropic and 
  $\lambda_\alpha + m_\alpha+2m_{2\alpha}$ is an even non-positive integer.
\item $\alpha$ is anisotropic and $\lambda_\alpha + m_\alpha+m_{2\alpha}$ is an odd negative integer.
\end{enumerate}
\end{CRL}

To set up the rank reduction, fix a positive system $\Sigma^+\subseteq\Sigma$. For any
positive system $\Phi\subseteq\Sigma$, we let
\[
\bar{\ger n}_\Phi\defi\bigoplus\nolimits_{\alpha\in(-\Sigma^+)\cap\Phi}\ger g^\alpha\ ,\ \ger n_\Phi\defi\bigoplus\nolimits_{\alpha\in\Sigma^+\cap\Phi}\ger g^\alpha\ ,\ \ger n^\Phi\defi\bigoplus\nolimits_{\alpha\in\Phi}\ger g^\alpha,
\]
and denote by $\bar N_\Phi$, $N_\Phi$, and $N^\Phi$ the respective analytic subsupergroups of $G$.

\begin{LMM}[nprod-decomp]
  Let $\Phi\subseteq\Sigma$ be a positive system. The embedding $\bar N_\Phi\to N^\Phi$ induces an isomorphism $\bar N_\Phi\cong N^\Phi/N_\Phi$.
\end{LMM}

\begin{PRF}
  By applying the inverse function theorem \cite{leites}, this follows from the equality $\ger n^\Phi=\ger n_\Phi\oplus\bar{\ger n}_\Phi$ and the classical case \cite{helgason84}*{Chapter IV, Lemma 6.8}.
\end{PRF}

\begin{CRL}[nbar-inv]
  The invariant Berezinian density $\Abs0{D\bar n}$ on $\bar N_\Phi$
  is invariant under the action of $N^\Phi$ induced by the isomorphism
  $\bar N_\Phi=N^\Phi/N_\Phi$.
\end{CRL}

\begin{PRF}
  On general grounds \cite{ah-berezin}*{Theorem 4.13}, \cite{a-hchom}*{Proposition A.2}, $\bar N_\Phi$ and $N^\Phi/N_\Phi$ have non-zero $\bar N_\Phi$- and $N^\Phi$-invariant Berezian densities $\Abs0{D\bar n}$ and $\Abs0{D\dot n}$, respectively, which are unique up to constant multiples. The pullback to $\bar N_\Phi$ of $\Abs0{D\dot n}$ is \emph{a fortiori} $\bar N_\Phi$-invariant, and hence, it differs from $D\bar n$ only by a non-zero scalar. The assertion follows.
\end{PRF}

\begin{DFN}
  For any $\alpha\in\Sigma$, we let $\alpha^+\defi\Sigma\cap(\rats_{>0}\cdot\alpha)$ and $\alpha^-\defi-\alpha^+$.
\end{DFN}

Let $\Phi$ and $\Psi$ be positive systems of $\Sigma$. For $\alpha\in B(\Phi)$ \scth
\[
\Psi=\alpha^-\cup\Phi\setminus\alpha^+,
\]
we write $\Phi\xrightarrow\alpha\Psi$ and say that $\Phi$ and $\Psi$ are \Define{adjacent}. For any positive system $\Phi\subseteq\Sigma$, there are positive systems $\Sigma^+=\Phi_0,\dotsc,\Phi_n=\Phi$ and $\alpha_i\in\Sigma^+\cap B(\Phi_i)$ \scth
\[
\Sigma^+\xrightarrow{\alpha_0}\dotsm\xrightarrow{\alpha_{n-1}}\Phi.
\]
If $n$ is minimal, then this is called a \Define{minimal gallery}.

Now, fix adjacent positive systems $\Phi\xrightarrow\alpha\Psi$ where $\alpha$ is indivisible and contained in $\Sigma^+$. Thus, $(-\Sigma^+)\cap \Psi=(-\Sigma^+)\cap\Phi\amalg\alpha^-$. Let $\bar N_\alpha$ be the analytic subsupergroup associated with the subalgebra $\bigoplus_{\gamma\in\alpha^-}\ger g^\gamma$ of $\bar{\ger n}$.

\begin{LMM}[alpha-mult]
  Multiplication induces an isomorphism $\bar N_\Psi\cong\bar N_\alpha\times\bar N_\Phi$ where $\bar N_\alpha$ is an analytic subsupergroup of $G$ with Lie superalgebra $\bar{\ger n}_\alpha=\bigoplus_{\beta\in\alpha^-}\ger g^\beta$.
\end{LMM}

\begin{PRF}
  This is proved in the same way as \REF{lmm}{nprod-decomp}.
\end{PRF}

\begin{CRL}[adj-nbar]
  Let $\Abs0{D\bar n}$, $\Abs0{D\bar n'}$, and $\Abs0{D\bar n''}$ denote the invariant Berezinian densities of $\bar N_\Psi$, $\bar N_\Phi$, and $\bar N_\alpha$, respectively. Then $\Abs0{D\bar n}=\Abs0{D\bar n''}\otimes\Abs0{D\bar n'}$.
\end{CRL}

\begin{PRF}
  The proof is analogous to that of \REF{crl}{nbar-inv}.
\end{PRF}

\begin{DFN}[weyl-ret]
  On $\bar N_\alpha$, we consider the standard Lie supergroup retraction, denoted by $r_\alpha$. We call retractions $r_\Phi$ on $\bar N_\Phi$ and $r_\Psi$ on $\bar N_\Psi$ \Define{related by $\alpha$} if 
  \[
    \phi_0\circ r_\Psi=(r_\alpha\times r_\Phi)\circ\phi, 
  \]
  where $\phi:\bar N_\Psi\to\bar N_\alpha\times\bar N_\Phi$ is the isomorphism from \REF{lmm}{alpha-mult}.

  Whenever we are given a minimal gallery 
  \[
    \Sigma^+=\Phi_0\xrightarrow{\alpha_0}\dotsm\xrightarrow{\alpha_{n-1}}\Phi_n=\Sigma^-,
  \]
  then the retraction on $\bar N=\bar N_{\Phi_n}$ obtained by the requirement that for any $i<n$, $r_{\Phi_i}$ be related by $\alpha_i$ to $r_{\Phi_{i+1}}$, is called a \Define{Weyl retraction}.
\end{DFN}

Consider the Weyl vectors
\[
\vrho_\Phi\defi\frac12\str_{\ger n^\Phi}\ad|_{\ger a}=\frac12\sum_{\beta\in\Phi}m_\beta\beta\nd 
\vrho_\alpha\defi\frac12\sum_{\beta\in\alpha^+}m_\beta\beta,
\]
where $m_\beta\defi\dim\ger g^\beta_\ev-\dim\ger g^\beta_\odd$. 

Assume that the retractions on $\bar N_\Phi$ and $\bar N_\Psi$ are related by $\alpha$. Whenever the integrals in question exist for $\lambda\in\ger a^*$, we set
\begin{equation}\label{eq:cphi}
  c_\Phi(\lambda)\defi\int_{\bar N_\Phi}e^{-(\lambda+\vrho_\Phi)(H(\bar n))}\,\Abs0{D\bar n}\ ,\ c_\alpha(\lambda)\defi\int_{\bar N_\alpha}e^{-(\lambda+\vrho_\alpha)(H(\bar n''))}\Abs0{D\bar n''},
\end{equation}
and similarly for $\Psi$.

\begin{PRO}[rankred]
  Retain the above assumptions. Then 
  \[
  c_\Psi(\lambda)=c_\Phi(\lambda)c_\alpha(\lambda)
  \]
  in the sense that the left-hand side converges absolutely if and only if the right-hand side does, and in this case, equality holds. Moreover, we have 
  \[
    c(\lambda+\vrho)=\int_{K/M}\frac{\Dual0{\mu_{-\lambda}}{\pi(k)v_\lambda}}{\Dual0{\mu_{-\lambda}}{v_\lambda}}\,\Abs0{D\dot k}.
  \]
\end{PRO}

The \emph{proof} is preceded by a number of technical statements.

\begin{PRO}[rho]
  Let $\Phi\xrightarrow\alpha\Psi$. Then $\Dual0{\vrho_\Psi}\alpha=\Dual0{\vrho_\alpha}\alpha$.
\end{PRO}

The \emph{proof} requires the following two lemmas.

\begin{LMM}[refl]
  Let $\alpha\in B(\Delta^+)$ be anisotropic and $\beta\in\Delta^+$ not proportional to $\alpha$. Then $r_\alpha(\beta)\in\Delta^+$. In particular, $r_\alpha(\Delta^+)=\Delta^+\setminus[\nats\alpha\cap\Delta]\cup[(-\nats\alpha)\cap\Delta]$.
\end{LMM}

\begin{PRF}
  Seeking a contradiction, we assume $r_\alpha(\beta)<0$. We have $r_\alpha(\beta)=\beta-n\alpha$ where $n=\smash{2\frac{\Dual0\alpha\beta}{\Dual0\alpha\alpha}}$ is an integer, and the root string $\beta+k\alpha$ has no gaps. (These standard facts may be derived by considering the sum of $\ger g_{\beta+k\alpha}$ as a module for the root algebra associated with $\alpha$, which is of type $\ger{sl}(2)$ or $\ger{osp}(1|2)$.) We may reduce to the situation that $\beta>0$ and $\beta-\alpha<0$. But then $\alpha=\beta+(\alpha-\beta)$ expresses $\alpha$ as the sum of two positive roots, contradiction.
\end{PRF}

\begin{LMM}[minrep]
  Let $\alpha\in\Phi$ and $\Delta^+_\Phi$ be an adapted positive system lying above $\Phi$. Any $\beta\in\Delta^+_\Phi$ with $\beta|_\ger a=\alpha$ is called a \Define{representative} of $\alpha$. Let $\sle$ be the order relation on $\Delta$ induced by $\Delta^+_\Phi$. If $\alpha$ is simple, then any minimal representative is simple.
\end{LMM}

\begin{PRF}
  Let $\alpha$ be simple and $\beta$ a minimal representative. Assume that $\beta=\gamma+\delta$, where $\gamma,\delta\in\Delta_\Phi^+$. Then $\alpha=\gamma|_\ger a+\delta|_\ger a$, where $\gamma|_\ger a,\delta|_\ger a\in\Phi\cup0$. Since $\alpha$ is simple, $\delta|_\ger a=0$ (w.l.o.g.). But this contradicts minimality of $\beta$, so $\beta$ is simple.
\end{PRF}

\begin{PRF}[\protect{Proof of \REF{pro}{rho}}]
  Let $\Delta^+_\Psi$ be an adapted positive system lying above $\Psi$. Then we have
  \begin{equation}\label{eq:rho-diff}
	2\Dual0{\vrho_\Psi-\vrho_\alpha}\alpha=\sum_{\beta\in\Delta^+_\Psi,\beta|_\ger a\notin\rats\alpha}(-1)^{\Abs0\beta}\Dual0\beta\alpha.
  \end{equation}
  Let $-\gamma\in\Delta_\Psi^+$ be a minimal representative of $-\alpha$, which is simple by \REF{lmm}{minrep}. By \REF{lmm}{refl} and the definition of odd reflections (see \REF{dfn}{oddrefl}), we find
  \[
  r_{-\gamma}(\Delta^+_\Psi)=\Delta^+_\Psi\setminus[(-\nats\gamma)\cap\Delta]\cup[(\nats\gamma)\cap\Delta].
  \]
  Let $\beta\in\Delta^+_\Psi$ \scth $\beta|_\ger a$ is not proportional to $\alpha$. Then $\beta$ is not proportional to $\gamma$, so $r_{-\gamma}(\beta)$ is in $\Delta^+_\Psi$.

  According to \cite{serganova-genroots}, the odd reflection $r_{-\gamma}(\Delta^+_\Psi)$ of $\Delta^+_\Psi$ is the image of $\Delta^+_\Psi$ under an involutive self-map of the set of roots, also denoted by $r_{-\gamma}$, given as follows:  If $\gamma$ is anisotropic, then $r_{-\gamma}(\beta)=r_{2\gamma}(\beta)$. This has the same parity as $\beta$, and
  \[
  \Dual0{r_{-\gamma}(\beta)}\gamma=\Dual0{\beta}{r_\gamma(\gamma)}=-\Dual0\beta\gamma.
  \]
  Moreover, we have
  \[
  r_\gamma(\beta)|_\ger a=\beta|_\ger a-n\alpha\notin\rats\alpha.
  \]

  If $\gamma$ is isotropic, then $r_{-\gamma}(\beta)=\beta\pm\gamma$. This root has parity opposite to $\beta$, and
  \[
  \Dual0{r_{-\gamma}(\beta)}\gamma=\Dual0{\beta\pm\gamma}\gamma=\Dual0\beta\gamma.
  \]
  Moreover, as before, $r_{-\gamma}(\beta)|_\ger a\notin\rats\alpha$. Thus, $r_{-\gamma}$ leaves the index set of the sum in Equation \eqref{eq:rho-diff} invariant, and
  \[
  	(-1)^{\Abs0{r_{-\gamma}(\beta)}}\Dual0{r_{-\gamma}(\beta)}\alpha=-(-1)^{\Abs0\beta}\Dual0\beta\alpha.
  \]
  Thereby, we find
  \[
  2\Dual0{\vrho_\Psi-\vrho_\alpha}\alpha=0,
  \]
  proving the proposition.
\end{PRF}

After these preliminaries, we prove the rank reduction formula.

\begin{PRF}[\protect{Proof of \REF{pro}{rankred}}]
  By \REF{crl}{adj-nbar} and \REF{crl}{nbar-inv}, and because of the equality $\log a(\bar n'')=H(\bar n'')$, we have
  \begin{align*}
    c_\Psi(\lambda)&=\int_{\bar N_\Phi\times\bar N_\alpha}e^{-(\lambda+\vrho)(H(a(\bar n'')\bar n'))}\Abs0{D\bar n'}\Abs0{D\bar n''}\\
    &=\int_{\bar N_\Phi\times\bar N_\alpha}e^{-(\lambda+\vrho)(H(a(\bar n'')\bar n'a(\bar n'')^{-1})}e^{-(\lambda+\vrho)(H(\bar n''))}\Abs0{D\bar n'}\Abs0{D\bar n''}\\
    &=\int_{\bar N_\Phi\times\bar N_\alpha}e^{-(\lambda+\vrho)(H(\bar n'))}e^{-(\lambda+\vrho+\str_{\bar{\ger n}_\Phi}\ad|_{\ger a})(H(\bar n''))}\Abs0{D\bar n'}\Abs0{D\bar n''}
    \intertext{%
      Next, we note that
      \begin{align*}
        2(\vrho+\str_{\bar{\ger n}_\Phi}\ad|_\ger a)&=\sum_{\beta\in\Sigma^+\cap\Phi}m_\beta\beta+\sum_{\beta\in(-\Sigma^+)\cap \Phi}m_\beta\beta
        =\sum_{\beta\in\Phi}m_\beta\beta=2\vrho_\Phi.
      \end{align*}
      Applying \REF{pro}{rho}, this gives
    }
    &=\int_{\bar N_\Phi}e^{-(\lambda+\vrho)(H(\bar n'))}\Abs0{D\bar n'}\int_{\bar N_\alpha}e^{-(\lambda+\vrho_\alpha)(H(\bar n''))}\Abs0{D\bar n''},
  \end{align*}
  which is the desired result.

  We now argue that the coordinate changes we have performed in transforming the integral from Equation \eqref{eq:kinv-nonzero} to the quantity $c(\lambda+\vrho)$ do not introduce any boundary terms. Indeed, taking the pullback along $k:\bar N\to K/M$ of some retraction of $K/M$, and then taking $\alpha$-related retractions in every step of the above recursion, we arrive by some retraction on $\bar N_\alpha$. Since we have pulled back retractions, no boundary terms can appear \cite{ahp-integration}. Moreover, if $\alpha$ is isotropic, $\bar N_\alpha$ is purely odd, so that the underlying space is compact, and we may change to the standard retraction on $\bar N_\alpha$ without introducing boundary terms.

  If $\alpha$ is anisotropic, then $\bar N_\alpha$ is the Iwasawa subgroup of some $\theta$-stable analytic subsupergroup $G_\alpha$, which satisfies the assumptions we have imposed on $G$. Thus, the integrand of $c_\alpha(\lambda+\vrho_\alpha)$ is the pullback along $k:\bar N_\alpha\to K_\alpha/M_\alpha$ of a Berezinian density similar to that considered in Equation \eqref{eq:kinv-nonzero}, only for the rank one pair $(G_\alpha,\theta)$.  The given retraction on $\bar N_\alpha$ extends to a retraction on $K_\alpha/M_\alpha$, since by construction, it is given by pullback from $K/M$. Since the underlying space of $K_\alpha/M_\alpha$ is again compact, it remains to show that passing from some retraction on $K_\alpha/M_\alpha$ to the standard retraction on $\bar N_\alpha$ does not introduce boundary terms. But this is done in \cite{ap-cfn}.
\end{PRF}

For the case of `rank one', we have the following result, \vq \cite{ap-cfn}.

\begin{PRO}[cfn-rk1]
  Assume that $\alpha$ is anisotropic. Identify $\lambda\in\ger a^*$ with $\lambda_\alpha\in\cplxs$. Then the integral $c_\alpha(\lambda)$ converges absolutely for $\Re\lambda>0$, and equals
  \[
  c_\alpha(\lambda)=c_0\frac{2^{-\lambda}\Gamma(\lambda)}{\Gamma\Parens1{\tfrac12(\tfrac{m_\alpha}2+1+\lambda)}\Gamma\Parens1{\tfrac12(\frac{m_\alpha}2+m_{2\alpha}+\lambda)}},
  \]
  for some non-zero constant $c_0$.
\end{PRO}

The case of an isotropic root is somewhat different. Indeed, if we apply the duplication formula for the $\Gamma$ function
to the right-hand side in \REF{pro}{cfn-rk1} for any $\alpha$, \scth $\Dual0\alpha\alpha=0$, then setting $m_\alpha=-2q <0$ yields
  \[
  \frac{2^{-\Dual0\lambda\alpha}\Gamma(\Dual0\lambda\alpha)}{\sqrt\pi\,2^{1+q-\Dual0\lambda\alpha}\Gamma\Parens1{\Dual0\lambda\alpha-q}}
  =\frac{2^{-q}}{2\sqrt\pi}\Parens1{\Dual0\lambda\alpha-1}\dotsm\Parens1{\Dual0\lambda\alpha-q}.
  \] 
Although this polynomial has the correct degree, the true value of $c_\alpha(\lambda)$ is different. 

\begin{PRO}[cfn-isotropic] Let $\alpha\in\Sigma^+$ be indivisible and
isotropic and $-2q:=m_\alpha$. Then for any $\lambda\in\ger a^*$ and a suitable
normalisation of Berezinians
  \[
c_\alpha(\lambda)=(-2)^{q}\Dual0\lambda\alpha^{q}.
  \]
\end{PRO}

\begin{PRF}
We have $2\alpha\notin\Sigma$, so $[\ger g^\alpha,\ger g^\alpha]=0$. Let
$\bar n\in_S\bar N_\alpha$, and write $\bar n=e^z$ where $z\in_S\ger
g^{-\alpha}$. To compute $H(\bar n)$, write $\bar n=kan$ and compute
\[
 \theta(\bar n)^{-1}\bar n=\theta(n)^{-1}a^2n .
\]
Writing  $a=e^{x h_\alpha}$ for $x\in_S\reals$, $b(h_\alpha,.)=\alpha$ and $n=e^y \in_S N_\alpha$ with $y\in_S\ger g^\alpha$, we have
\[ 
\theta(n)^{-1}a^2n=e^{-\theta(y)}e^ye^{2xA_\alpha}
\]
since $[h_\alpha, y]=\Dual0\alpha\alpha\cdot y=0$. Further $\mf g_\alpha
\subseteq \mf g_\odd$, hence we have $y^2=0=z^2$ in $\Gamma(\sh O_S)\otimes\Uenv0{\ger
g}$.
 Moreover, $[[\ger g^\alpha,\ger g^{-\alpha}],\ger
g^{\pm\alpha}]=0$ by \cite{a-hchom}*{Lemma 3.5}. Hence, the
Campbell--Hausdorff formula on $S$-valued points gives
\[
\theta(n)^{-1}a^2n=e^{-\theta(y)+y-[\theta(y),y]}e^{2xA_\alpha}=e^{-\theta(y)+y-[\theta(y),y]+2xA_\alpha}
\] and similarly
\[ \theta(\bar n)^{-1}\bar n=e^{-\theta(z)+z-[\theta(z),z]}.
\] Comparing these quantities, we find $y=-\theta(z)$ and hence
\[ -[\theta(z),z]=-[\theta(y),y]+2xA_\alpha=-[z,\theta(z)]+2xA_\alpha,
\] so that
\[ H(\bar n)=[z,\theta(z)]_\ger p=[z,\theta(z)],
\] where we denote by $(-)_\ger p$ the $\ger p$-projection
$\tfrac12(\id-\,\theta)$.

We now introduce coordinates on $\bar N_\alpha$. Observe that
$b^\theta(x,y)=b(x,\theta y)$ induces a symplectic form on $\ger
g^\alpha$. Let $2q\defi-m_\alpha$ and choose a symplectic basis
$\xi_i,\eta_i$ of $\ger g^\alpha$, so that
\[
 b^\theta(\xi_i,\eta_j)=\delta_{ij}\ ,\
b^\theta(\xi_i,\xi_j)=b^\theta(\eta_i,\eta_j)=0.
\]
 We define $\bar\xi_i\defi\theta(\xi_i)$ and
$\bar\eta_i\defi\theta(\eta_i)$. This gives a symplectic basis of
$\ger g^{-\alpha}$.

From \cite{a-hchom}*{Lemma 3.5}, we deduce the bracket relations
\[ [\xi_i,\bar\eta_j]_\ger p=-[\eta_i,\bar\xi_j]_\ger
p=\delta_{ij}A_\alpha\ ,\ [\xi_i,\bar\xi_j]_\ger
p=[\eta_i,\bar\eta_j]_\ger p=0.
\] Writing $\bar n=e^z$, where
\[ z=\sum_ia^i\bar\xi_i+b^i\bar\eta_i
\] with $a^i,b^i\in\Gamma(\sh O_S)_\odd$, we find
\begin{align*} H(\bar
n)&=\sum_{ij}a^ia^j[\xi_i,\bar\xi_j]+a^ib^j[\eta_i,\bar\xi_j]+b^ia^j[\xi_i,\bar\eta_j]+b^ib^j[\eta_i,\bar\eta_j]\\
&=\sum_i(a^ib^i-b^ia^i)A_\alpha=2\sum_ia^ib^iA_\alpha.
\end{align*} Since $\Dual0{\vrho_\alpha}\alpha=0$, we obtain after
suitable normalisation
\[ c_\alpha(\lambda)=\Parens3{\int
D(\xi,\eta)\,(1-2\Dual0\lambda\alpha\xi\eta)}^q=(-2)^q\Dual0\lambda\alpha^q,
\]
which proves the claim.
\end{PRF}

\begin{RM}
  Observe that for $\Dual0\alpha\alpha=0$, we have $m_\alpha=-2$ in all known examples, so the zeros in \REF{pro}{cfn-isotropic} are still simple. 
\end{RM}

\section{Self-dual $\mf{gl}$ representations}
\label{sec:selfdual}

In this section, as an application of \REF{thm}{main}, we discuss the self-duality of highest weight representations for the special case of a symmetric superpair of $\ger{gl}$ type. Throughout this section, let $p,q,r,s \in \N$ such that
\[
  (p-q)(r-s)\geq0.
\]
As mentioned in \REF{eg}{glnn}, this ensures that $\mf g = \mf
{gl}^{p+q|r+s}$ is of even type, and hence by \REF{eg}{gl-iwasawa}, 
we can apply \REF{thm}{main}.

W.l.o.g., we may assume that $r\geq s$ and $p\geq q$, as will become clear in \REF{dfn}{glpqrs}.
We consider $\g := \gl^{p+q|r+s}\simeq \uEnd(\C^{p+q|r+s})$ with
standard basis $E_{i,j} := e_i\otimes e_j^*$.
Before we cast
the definitions from \REF{eg}{glnn} in an explicit matrix from in
\REF{sec}{compSys}, let us start by reviewing the concept of
$\delta\eps$-chains from \cite{cheng_wang}.  These will give us a parametrisation of positive systems, which is more convenient than that given by the Dynkin diagrams.

\subsection{$\delta \eps$-chains}

Our main tool in discussing self-duality will be a chain of simple reflections mapping
the simple system $\Pi$ to $- \Pi$. We will describe this chain
explicitly in this subsection. It is similar to the corresponding chain of ordinary reflections, only that some of the simple reflections are now odd, see \REF{sec}{oddRefl}. To give a full description, we introduce $\delta\eps$-chains.

\begin{DFN}
  Consider a string $S_1 \ldots S_{n}$ with $S_k \in \{\delta_i,
  \eps_j\}$ and $S_i \neq S_j$ for $i\neq j$. Any such string is
  called a \emph{$\delta\eps$-chain}. 
We will use the notation
  \[
  \Pi(S_1 \ldots S_{n}) \defi  \{S_1 - S_2, S_2 - S_3, \ldots , S_{n-1} - S_{n}\}.
  \]
\end{DFN}

\begin{LMM}[simple-swapp]
  If $\Pi(S_1 \ldots S_{n})$ is a system of simple roots, then simple
  reflections correspond to swapping neighbouring elements of the $\delta\eps$-chain, \ie
  \[
  r_{S_i - S_{i+1}} \left(\Pi\left(S_1 \ldots S_{n}\right)\right) =  \Pi(S_1 \ldots
  S_{i-1}S_{i+1}S_{i}S_{i+2} \ldots S_{n}).
  \]
\end{LMM}

\begin{PRF}
  If $r_{S_i - S_{i+1}}$ is an even reflection then
  \begin{align}
    r_{S_i - S_{i+1}}(S_k) &= S_k - 2 \frac{\Dual0{S_i - S_{i+1}}{S_k}}{\Dual0{S_i -
      S_{i+1}}{S_i - S_{i+1}}} (S_i - S_{i+1})
    \\
    &=
    S_k - (\delta_{k,i} - \delta_{k, i+1})(S_i - S_{i+1})
    =
    \begin{cases}
      S_k & k \notin \{i,i+1\},\\
      S_{i+1} & k=i,\\
      S_i & k=i+1.
    \end{cases}
  \end{align}
  %
  If $r_{S_i - S_{i+1}}$ is odd, then $r_{S_i - S_{i+1}}$ replaces
  \[
  (S_k- S_{k+1})
  \mapsto
  \begin{cases}
    S_k - S_{k+1} & k \notin \{i-1,i,i+i\},\\
    S_k - S_{k+1} + S_i - S_{i+1}= S_{i-1} - S_{i+1} & k=i-1,\\
    S_k - S_{k+1} + S_i - S_{i+1}= S_i - S_{i+2} & k=i+1,\\
    S_{i+1} - S_{i} & k =i,
  \end{cases}
  \]
by \REF{dfn}{oddReflOfHW}.
\end{PRF}

\begin{LMM}
Let $\delta_i,\eps_j\in \mf h ^*$ as defined in Equation \eqref{eq:defDeltaEpsilon}.
  Then the systems of simple roots, $\Pi(C)$, of $\gl^{p+q|r+s}$ are in one-to-one
  correspondence with the $\delta \eps$-chains $C$ of full length
  $p+q+r+s$, \ie containing all $\delta_i$ and $\eps_j$.
\end{LMM}

\begin{PRF}
  As shown in \REF{lmm}{posSystem}, Equation \eqref{eq:compatibleChain}, there is a
  $\delta\eps$-chain corresponding to a specific system of simple roots.
  By \cite{cheng_wang}*{Corollary 1.27} any fundamental system $\Pi'$
  can be produced by applying a chain of simple reflections to $\Pi$ and
  any such yields a fundamental system.
As we have seen in \REF{lmm}{simple-swapp}, simple reflections 
amount to swapping neighbours in a $\delta\eps$-chain. These swaps generate
all permutations and hence all chains.
\end{PRF}

\begin{CRL}
  If $\Pi$ corresponds to a certain $\delta\eps$-chain, then $-\Pi$
  corresponds to the reversed chain.
\end{CRL}

\begin{LMM}[BubbleSort]
  A chain of simple reflections reverting a $\delta\eps$-chain $C=S_1  S_2 S_3 \ldots  S_n$, that is, mapping $\Pi(C)$ to $-\Pi(C)$, is the following:
  \begin{align}
    R_C \defi \ldots &
    \circ
    \left( r_{S_{3}-S_{n-1}} \circ\ldots\circ r_{S_{n-3}-S_{n-1}} \circ r_{ S_{n-2}-S_{n-1}} \right)
    \\&\circ
    \left( r_{S_2 - S_{n-1}} \circ\ldots\circ r_{S_2 - S_4} \circ r_{S_2 - S_3} \right)
    \\&\circ
    \left(r_{S_2}- r_{S_{n}} \circ\ldots\circ r_{S_{n-2}-S_n} \circ r_{S_{n-1}-S_n} \right)
    \\&\circ
    \left( r_{S_1 - S_n} \circ\ldots\circ r_{S_1 - S_3} \circ r_{S_1 - S_2} \right)
    .\textbf{}
  \end{align}
\end{LMM}


\subsection{Restricted roots for
  $\mf{gl}^{p+q|r+s}$}
\label{sec:compSys}

We fix an involution on $\ger g$ and compute the corresponding restricted root data.

\begin{DFN}[glpqrs]
  The involution defining $\g = \mf{k}\oplus \mf{p}$ is $\theta: X
  \mapsto  \sigma X \sigma$ with
  \begin{align}
    \sigma &=
    \left(
      \begin{array}{cc|cc}
        \1_p &0 & 0 & 0\\
        0&-\1_q & 0 & 0 \\
        \hline
        0 & 0      &\1_r &0\\
        0 & 0      &      0&-\1_s
      \end{array}
    \right).
    &\text{Hence }
    \g&=
    \left\{
      \left(
        \begin{array}{cc|cc}
          \mf{k}_\ev & \p_\ev & \mf{k}_\odd & \p_\odd\\
          \p_\ev & \mf{k}_\ev & \p_\odd & \mf{k}_\odd\\
          \hline
          \mf{k}_\odd & \p_\odd & \mf{k}_\ev & \p_\ev\\
          \p_\odd &\mf{k}_\odd & \p_\ev & \mf{k}_\ev
        \end{array}
      \right)
    \right\}
    \label{glIntoHP}
  \end{align}
  is the $\theta$-decomposition of $\mf g$, \ie $\mf{k} =\mathfrak{gl}^{p|r} \oplus\gl^{q|s}$.
  The even $\uEnd(\C^{p+q})$ and $\uEnd(\C^{r+s})$ parts will be called the \emph{boson-boson} and \emph{fermion-fermion} blocks, respectively.
\end{DFN}

\begin{LMM}[NonCompactCSFormOfgl]
  A non-compact cs form of $(\g,\theta)$ is given by
  \[
  \g_{\ev,\R}=
  \left\{
    \left(
      \begin{array}{cc|cc}
        A & B & 0 & 0 \\
        B^\dagger & C & 0 & 0 \\
        \hline
        0&0 & D &E\\
        0&0 & E^\dagger &F
      \end{array}
    \right)
    =
    \left(
      \begin{array}{cc|cc}
        -A^\dagger & B & 0 & 0 \\
        B^\dagger & -C^\dagger & 0 & 0 \\
        \hline
        0&0 & -D^\dagger &E\\
        0 &0 & E^\dagger & -F^\dagger
      \end{array}
    \right)
  \right\}
  \]
\end{LMM}

\begin{PRF}


  This corresponds to the non-compact Lie algebras $\mf u(p,q)$ and $\mf u(r,s)$.
\end{PRF}

\begin{LMM}[glEvenType]
  We can choose a real even Cartan subspace subspace $\mfa_{\ev,\R}
  \subseteq\p_{\ev,\R}$ using the following notation
  \begin{align}
    &A :
    \C^q\oplus\C^s \to \uEnd(\C^{p+q})\oplus\uEnd(\C^{r+s}) \hookrightarrow \g
    \\
    &A(a^B+a^F) :=
    \left(
      \begin{array}{cc|cc}
        0& \diag(a^B)&0&0\\
        \diag(a^B)^\dagger& 0 & 0 & 0 \\
        \hline
        0&0 &  0 & \diag(a^F)\\
        0 &0&  \diag(a^F)^\dagger &0
      \end{array}
    \right)
  \end{align}
  to set $\mfa_{\ev,\R} := A (i \R^q \oplus i \R^s)$.
\end{LMM}

\begin{LMM}[csa-choice]
  The following choice of $\mf{h}$ defines an even Cartan subalgebra of $\g$ with
  $\mfa \subseteq \mathfrak{h}$. We use
  \begin{align}
    &m :
    \C^q \oplus \C^{p-q}\oplus\C^s\oplus\C^{r-s} \to \mf k
    \\
    &m(b^B+c^B+b^F+c^F) :=
    \diag\left(b^B,c^B,b^B|b^F,c^F,b^F\right)
  \end{align}
  to specify $\mf h_{0,\R} := m(i\R^q \oplus i\R^{p-q} \oplus i \R^{r} \oplus i \R^{r-s})
  \oplus \mf a_{0,\R}$.
  This also defines linear coordinates on $\mf{h}$ with respect to a
  super trace orthonormal basis.
\end{LMM}

\subsubsection{Roots}
Diagonalising the adjoint action of $\mf h$ on $\g$ is a trivial
computational task. We record the results as follows.

\begin{LMM}[BBroots]
  The even roots in the boson-boson-block and corresponding root vectors
  are given by

  \begin{align}
    \begin{array}{c|c|c}
      \text{root } \alpha & \g^\alpha \text{ spanned by} & \text{for}\\
      \hline
      c^{B}_{i-q} - b^{B}_{j} \pm i  a^B_{j}
      & E_{i,j} \pm i E_{i,j+p}
      & q < i \leq p, j\leq q \\
      b^{B}_{i} - c^{B}_{j-q} \pm i a_{i}^B
      & E_{i,j} \pm i E_{i+p,j}
      &  i \leq q,  q<j\leq p \\
      b^{B}_{i} - b^{B}_{j} \pm i \left(a_{i}^B - a_{j}^B\right)
      &
      \begin{array}{c}
        E_{i,j} + E_{i+p,j+p}\\
        \pm i \left(E_{i+p,j}-E_{i,j+p}\right)
      \end{array}
      &  i\neq j\leq q \\
      b^{B}_{i} - b^{B}_{j} \pm i \left(a_{i}^B + a^B_{j}\right)
      &
      \begin{array}{c}
        E_{i,j} - E_{i+p,j+p} \\
        \pm i \left(E_{i+p,j} +E_{i,j+p}\right)
      \end{array}
      & i,j \leq q \\
      c^{B}_{i-q} - c^{B}_{j-q}
      & E_{i,j}
      & q < i \neq j \leq p
    \end{array}
  \end{align}
  Note that in particular $\pm 2 i a^B_i$ is a root for $q>0$.

\end{LMM}

\begin{LMM}
  Applying the correspondence $B\mapsto F$, $q\mapsto s$, and $p \mapsto r$,
  to \REF{lmm}{BBroots}, yields the even roots in the fermion-fermion-block.
\end{LMM}

\begin{LMM}
  The odd roots with root vectors in the boson-fermion-block are the following.
  \begin{align}
    \begin{array}{c|c|c}
      \text{root } \alpha & \g^\alpha \text{ spanned by} & \text{for}\\
      \hline
      b_i^B \pm i (a_i^B - a_j^F) - b_j^F
      &
      \begin{array}{c}
        E_{i,j+p+q}  + E_{i+p,j+p+q+r}\\
        \pm i (E_{i+p,j+p+q} - E_{i,j+p+q+r})
      \end{array}
      & i\leq q, j\leq s
      \\
      b_i^B \pm i (a_i^B + a_j^F) - b_j^F
      &
      \begin{array}{c}
        E_{i,j+p+q} - E_{i+p,j+p+q+r}\\
        \pm i (E_{i+p,j+p+q} - E_{i,j+p+q+r})
      \end{array}
      & i\leq q, j\leq s
      \\
      b_i^B - c_{j-s}^F \pm i a_i^B
      & E_{i,j+p+q} \pm i E_{i+p,j+p+q}
      & i \leq q, s<j\leq r
      \\
      c_{i-q}^B - b_j^F \pm i a_j^F
      & E_{i,j+p+q} \pm i E_{i,j+p+q+r}
      & q<i \leq p, j\leq s
      \\
      c_{i-q}^B - c_{j-s}^F
      & E_{i,j+p+q}
      & q<i \leq p, s<j\leq r
    \end{array}
  \end{align}

  The remaining odd roots with root vectors in the fermion-boson-block are the following.
  \begin{align}
    \begin{array}{c|c|c}
      \text{root } \alpha & \g^\alpha \text{ spanned by} & \text{for}\\
      \hline
      b_i^F \pm i (a_i^F - a_j^B) - b_j^B
      &
      \begin{array}{c}
        E_{i+p+q,j} + E_{i+p+q+r,j+p}\\
        \pm i (E_{i+p+q+r,j} - E_{i+p+q,j+p})
      \end{array}
      & i\leq s, j\leq q
      \\
      b_i^F \pm i (a_i^F + a_j^B) - b_j^B
      &
      \begin{array}{c}
        E_{i+p+q,j}  - E_{i+p+q+r,j+p}\\
        \pm i (E_{i+p+q+r,j} - E_{i+p+q,j+p})
      \end{array}
      & i\leq s, j\leq q
      \\
      b_i^F - c_{j-q}^B \pm i a_i^F
      & E_{i+p+q,j} \pm i E_{i+p+q+r,j}
      & i \leq s, q<j\leq p
      \\
      c_{i-s}^F - b_j^B \pm i a_j^B
      & E_{i+p+q,j} \pm i E_{i+p+q,j+p}
      & s<i \leq r, j\leq q
      \\
      c_{i-s}^F - c_{j-q}^B
      & E_{i+p+q,j}
      & s<i \leq r, q<j\leq p
    \end{array}
  \end{align}

\end{LMM}

\subsubsection{Compatible positive system}

We fix a positive system of $\ger g=\ger{gl}^{p+q|r+s}$, compatible with our choice of $\ger a$. 

Let $\Pi$ be simple system corresponding to the $\delta\eps$-chain
  \begin{equation}
    \label{eq:compatibleChain}
    \delta_{p+1}\ldots\delta_{p+q}
    \eps_{r+1}\ldots\eps_{r+s}
    \delta_{q+1}\ldots\delta_{p}
    \eps_{s+1}\ldots\eps_r
    \eps_{s}\ldots\eps_{1}
    \delta_{q}\ldots\delta_{1}
  \end{equation}
  where
  \begin{align}
    \label{eq:defDeltaEpsilon}
    \begin{array}{ll|ll}
      \delta_j \defi b_j^B - i a_j^B & j\leq q
      &\eps_j \defi b^F_j-i a^F_j & j\leq s\\
      \delta_{j} \defi c_{j-q}^B & q<j\leq p
      &\eps_{j} \defi c^F_{j-s} & s<j\leq r\\
      \delta_{p+j} \defi b_j^B + i a_j^B & j\leq q
      &\eps_{r+j} \defi b_j^F+i a_j^F & j\leq s
    \end{array}
  \end{align}

\begin{LMM}[posSystem]
  The positive system $\Delta^+$ defined by $\Pi$ is compatible with $\theta$.
\end{LMM}

\begin{proof}
  We have that 
  \begin{align}
    \label{eq:SigmaPlus}
    \Delta^{+}|_{\mf{a}} \setminus \{0\} &=\Sigma^{+}(\mf{g}:\mf{a})
    \\
    &=
    \begin{aligned}[t]
      &\Braces1{i a^B_k, i(a_k^B +a_\ell^B), i a^F_k, i(a_k^F+a_\ell^F),
        i (a^B_k - a^F_\ell),i (a^B_k + a^F_\ell)} \\
      &\cup
      \Set1{i (a^B_k - a^B_\ell), i (a^F_k - a^F_\ell)}{k<\ell},      
    \end{aligned}
  \end{align}
  where it is understood that the short roots $ia^B_k$ and $ia^F_\ell$ are present only in the case that $q<p$ or $s<r$. This defines a positive system of the system $\Sigma$ of restricted roots, proving the claim. 
\end{proof}


  %
  In case $q<p$ and $s<r$, $\Pi$ corresponds to the Dynkin diagram
  
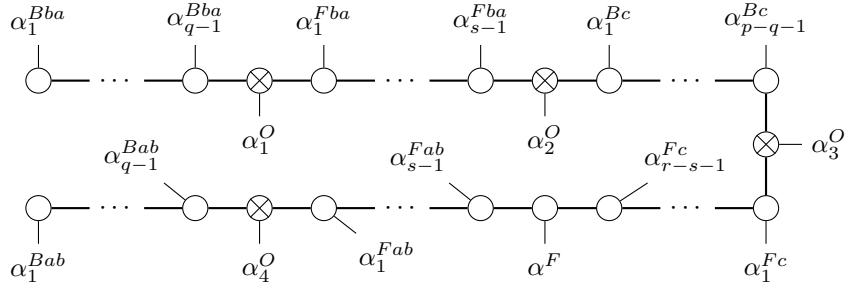
\begin{center}
    \begin{tikzpicture}
      \pgfmathsetlengthmacro{\dia}{0.2cm};
      \pgfmathsetlengthmacro{\dist}{0.5cm};
      \pgfmathsetlengthmacro{\Ldist}{0.3cm};
      \node[draw, circle, minimum width=\dia](N1){};
      \node[above=\Ldist of N1] {$\alpha_1^{Bba}$} edge (N1);

      \node[right = \dist of N1] (D1){$\cdots$} edge[thick] (N1);

      \node[draw, circle, minimum width=\dia, right = \dist of D1](Ns-1){} edge[thick] (D1);
      \node[above=\Ldist of Ns-1] {$\alpha_{q-1}^{Bba}$} edge (Ns-1);

      \node[draw, circle, minimum width=\dia, right = \dist of Ns-1](Ns){} edge[thick] (Ns-1);
      \draw (Ns.north east)--(Ns.south west);
      \draw (Ns.north west)--(Ns.south east);
      \node[below=\Ldist of Ns] {$\alpha_1^{O}$} edge (Ns);

      \node[draw, circle, minimum width=\dia, right = \dist of Ns](Ns+1){} edge[thick] (Ns);
      \node[above=\Ldist of Ns+1] {$\alpha_{1}^{Fba}$} edge (Ns+1);

      \node[right = \dist of Ns+1] (D2){$\cdots$} edge[thick] (Ns+1);

      \node[draw, circle, minimum width=\dia, right = \dist of D2](Ns+q-1){} edge[thick] (D2);
      \node[above=\Ldist of Ns+q-1] {$\alpha_{s-1}^{Fba}$} edge (Ns+q-1);

      \node[draw, circle, minimum width=\dia, right = \dist of Ns+q-1](Ns+q){} edge[thick] (Ns+q-1);
      \draw (Ns+q.north east)--(Ns+q.south west);
      \draw (Ns+q.north west)--(Ns+q.south east);
      \node[below=\Ldist of Ns+q] {$\alpha_{2}^{O}$} edge (Ns+q);

      \node[draw, circle, minimum width=\dia, right = \dist of Ns+q](Ns+q+1){} edge[thick] (Ns+q);
      \node[above=\Ldist of Ns+q+1] {$\alpha_{1}^{Bc}$} edge (Ns+q+1);

      \node[right = \dist of Ns+q+1] (D3){$\cdots$} edge[thick] (Ns+q+1);

      \node[draw, circle, minimum width=\dia, right = \dist of D3](Nr+q-1){} edge[thick] (D3);
      \node[above=\Ldist of Nr+q-1] {$\alpha_{p-q-1}^{Bc}$} edge (Nr+q-1);

      \node[draw, circle, minimum width=\dia, below = \dist of Nr+q-1](Nr+q){} edge[thick] (Nr+q-1);
      \draw (Nr+q.north east)--(Nr+q.south west);
      \draw (Nr+q.north west)--(Nr+q.south east);
      \node[right=\Ldist of Nr+q] {$\alpha^{O}_3$} edge (Nr+q);

      \node[draw, circle, minimum width=\dia, below = \dist of Nr+q](Nr+q+1){} edge[thick] (Nr+q);
      \node[below=\Ldist of Nr+q+1] {$\alpha_{1}^{Fc}$} edge (Nr+q+1);

      \node[left = \dist of Nr+q+1] (D4){$\cdots$} edge[thick] (Nr+q+1);

      \node[draw, circle, minimum width=\dia, left = \dist of D4](Nr+q+s-1){} edge[thick] (D4);
      \node[above right=\Ldist of Nr+q+s-1] {$\alpha^{Fc}_{r-s-1}$} edge (Nr+q+s-1);

      \node[draw, circle, minimum width=\dia, left = \dist of Nr+q+s-1](M1){} edge[thick] (Nr+q+s-1);
      \node[below=\Ldist of M1] {$\alpha^{F}$} edge (M1);

      \node[draw, circle, minimum width=\dia, left = \dist of M1](M2){} edge[thick] (M1);
      \node[above left=\Ldist of M2] {$\alpha^{Fab}_{s-1}$} edge (M2);

      \node[left = \dist of M2] (D6){$\cdots$} edge[thick] (M2);

      \node[draw, circle, minimum width=\dia, left = \dist of D6](M3){} edge[thick] (D6);
      \node[below right=\Ldist of M3] {$\alpha^{Fab}_{1}$} edge (M3);

      \node[draw, circle, minimum width=\dia, left = \dist of M3](M4){} edge[thick] (M3);
      \draw (M4.north east)--(M4.south west);
      \draw (M4.north west)--(M4.south east);
      \node[below=\Ldist of M4] {$\alpha^{O}_{4}$} edge (M4);

      \node[draw, circle, minimum width=\dia, left = \dist of M4](M5){} edge[thick] (M4);
      \node[above left=\Ldist of M5] {$\alpha^{Bab}_{q-1}$} edge (M5);

      \node[left = \dist of M5] (D7){$\cdots$} edge[thick] (M5);

      \node[draw, circle, minimum width=\dia, left = \dist of D7](M6){} edge[thick] (D7);
      \node[below=\Ldist of M6] {$\alpha^{Bab}_{1}$} edge (M6);

    \end{tikzpicture}
  \end{center}


  From this, the Dynkin diagram for restricted simple system is obtained by removing the
  $a_1^{Bc}$ to $a^{FC}_{r-s-1}$ part and `folding' (identifying 
  vertical pairs):

  \smallskip
  \begin{center}
    \begin{tikzpicture}
      \pgfmathsetlengthmacro{\dia}{0.2cm};
      \pgfmathsetlengthmacro{\dist}{0.5cm};
      \pgfmathsetlengthmacro{\Ldist}{0.3cm};
      \node[draw, circle, minimum width=\dia](N1){};

      \node[right = \dist of N1] (D1){$\cdots$} edge[thick] (N1);

      \node[draw, circle, minimum width=\dia, right = \dist of D1](Ns-1){} edge[thick] (D1);

      \node[draw, circle, minimum width=\dia, right = \dist of Ns-1](Ns){} edge[thick] (Ns-1);
      \draw (Ns.north east)--(Ns.south west);
      \draw (Ns.north west)--(Ns.south east);

      \node[draw, circle, minimum width=\dia, right = \dist of Ns](Ns+1){} edge[thick] (Ns);

      \node[right = \dist of Ns+1] (D2){$\cdots$} edge[thick] (Ns+1);

      \node[draw, circle, minimum width=\dia, right = \dist of D2](Ns+q-1){} edge[thick] (D2);

      \node[fill, draw, circle, minimum width=\dia, right = \dist of Ns+q-1](Ns+q){};
      \draw[thick] (Ns+q-1.north east) -- (Ns+q.north west);
      \draw[thick] (Ns+q-1.south east) -- (Ns+q.south west);
      \draw[thick] (Ns+q.west) -- +(-\dia,\dia);
      \draw[thick] (Ns+q.west) -- +(-\dia,-\dia);

    \end{tikzpicture}
  \end{center}


\subsection{Self-duality} 

We can finally characterise self-duality of highest weight modules. To that end, we introduce a decomposition of the highest weights $\lambda$, which is adapted to the choice of positive system made above. 

  To that end, note that $\{\delta_i,\eps_j\}$, as defined in Equation \eqref{eq:defDeltaEpsilon},
  form a basis of $\mf h^*$, where $\ger h\ni\1$ is the Cartan algebra of $\ger g=\ger{gl}^{p+q|r+s}$ defined in \REF{lmm}{csa-choice}. It is orthonormal in
  the sense that $\Dual0{\delta_i}{\delta_j} =  \delta_{i,j}$,
  $\Dual0{\delta_i}{\eps_j} = 0$ and
  $\Dual0{\eps_i}{\eps_j} = -  \delta_{i,j}$.

  For later reference, we give an explicit parametrisation of the highest weights. To that end, we order the orthonormal basis according to Equation \eqref{eq:compatibleChain} and write 
  \[
    \left(
      \lambda^{\delta}_{p+1},\ldots,\lambda^{\delta}_{p+q}
      \big |
      \lambda^{\eps}_{r+1},\ldots,\lambda^\eps_{r+s}
      \big |
      \lambda^{\delta}_{q+1},\ldots,\lambda^{\delta}_{p}
      \big |
      \lambda^{\eps}_{s+1},\ldots,\lambda^\eps_{s+r}
      \lambda^{\eps}_{s},\ldots,\lambda^\eps_{1}
      \big |
      \lambda^{\delta}_{q},\ldots,\lambda^{\delta}_{1}
    \right)
  \]
to represent
  \[
    \lambda=\sum_i \lambda^{\delta}_i \delta_i + \sum_j \lambda^{\eps}_j \eps_j.
  \]
  A glance at Equation \eqref{eq:defDeltaEpsilon} reveals that the $\lambda \in \mf{a}^*$ have the form
  \begin{equation}\label{eq:lambda-astar}
  \lambda=\left(
    -\lambda^\delta_1,\ldots,-\lambda_q^\delta
    \big |
    -\lambda^\eps_1,\ldots,-\lambda_s^\eps
    \big | 0\ldots0 \big | 0\ldots0
    \lambda^\eps_s, \ldots, \lambda_1^\eps
    \big |
    \lambda^\delta_q, \ldots, \lambda_1^\delta
  \right).
  \end{equation}

The following lemma makes the conditions stated in \REF{crl}{conditionsOnHW} explicit.

\begin{LMM}[sphericalHeighestWeights]
  Let $\lambda\in\ger a^*$ be the highest weight of a finite-dimensional irreducible spherical
  $\gl^{p+q|r+s}$ representation. With coefficients defined as above, we have that 
  \[
  -\lambda^\delta_1
  \geq \ldots\geq
  -\lambda^\delta_q
  \geq
  \lambda_1^\eps \geq \lambda^\eps_{2} \geq \ldots \geq
  \lambda^\eps_s \geq 0.
  \]
  Moreover, if $q<p$ or $s<r$, then these numbers are even, \ie $\lambda^S_{j} \in 2 \Z$. Otherwise, they are either all even or all odd. 
\end{LMM}

\begin{PRF}
  We have $i a_i^B = \frac{1}{2}(\delta_{p+i} - \delta_i)$ and $i a_j^F =
  \frac{1}{2}(\eps_{r+j}-\eps_j)$, hence
  \[
  \mf a^* = \Braces2{
    \textstyle\sum_{i=1}^q \lambda^{\delta}_i (\delta_i-\delta_{p+i}) + \sum_{j=1}^s
    \lambda^{\eps}_j (\eps_j - \eps_{r+j})
  }.
  \]
  For spherical representations, we have that $\forall \alpha_i \in
  \Sigma^+ : \Dual0{\lambda}{ \alpha_i} \in 2\N$, with $\Sigma^+$ given in
  Equation \eqref{eq:SigmaPlus}.
  Hence $\Dual0{i a_i^B}{\lambda} = - \lambda_i^\delta \in 2 \N$ and
  $\Dual0{i a_j^F}{\lambda} = \lambda_j^\eps \in 2 \N$ where we recall
  $\Dual0{\eps_i}{\eps_i} = -1$ which leads to the sign difference.
  Further
  \[
  \Dual0{i (a_i^B-a_{i+1}^B)}{\lambda} = \lambda_{i+1}^\delta -
  \lambda_i^\delta \in 2 \N
  \Rightarrow \lambda_{i+1}^\delta \geq \lambda_i^\delta
  \]
  and similarly with reversed signs for fermions. Finally
  \[
  \Dual0{i (a_i^B-a_{j}^F)}{\lambda} = -\lambda_{i}^\delta -
  \lambda_j^\eps \in 2 \N
  \Rightarrow -\lambda_{j}^\eps \geq \lambda_i^\delta
  .
  \]
  This proves our claim.
\end{PRF}


We now come to the main result of this section.

\begin{DFN}
  Let $R$ be a chain of even and odd reflections. Then we denote by $R'$ the chain
  which is obtain from $R$ by dropping all odd reflections.

  We call a $\delta\eps$-chain a
  \emph{palindrome} if the chain of Greek letters upon dropping the
  indices becomes a palindrome.
\end{DFN}

\begin{PRO}[mainSelfDual]
  Let $P$ be a $\delta\eps$-chain and $R_P$ the corresponding
  chain of reflections defined in \REF{lmm}{BubbleSort}.
  Then $R_P = R_P'$ as operators on $\mf h^*$ if and only if $P$
  is a palindrome.
\end{PRO}


\begin{PRF}
  Let $P$ be a palindrome.
  Without loss of generality, $P$ starts with $\delta$.
  By definition $r_{\eps-\delta} \circ r_{\delta-\eps} = \id$,
  which immediately leads to
  \begin{equation}
    \label{eq:elementPair}
    r_{\eps-\delta_1}\circ r_{\delta_{1}-\delta_{2}} \circ \ldots
    \circ r_{\delta_{k-1}-\delta_{k}}\circ r_{\delta_k-\eps} = r_{\delta_{1}-\delta_{2}} \circ \ldots
    \circ r_{\delta_{k-1}-\delta_{k}}.
  \end{equation}
  This already settles the case of $P=\delta \eps_1 \ldots
  \eps_n \delta'$.

  Now, we proceed inductively with $P=\delta Q\delta'$, where
  $Q=s_1\ldots s_n$ is any palindrome. Denote
  \[
    R(\delta)=r_{\delta-\delta'}\circ r_{\delta-s_n}\circ\ldots\circ r_{\delta-s_1}\nd 
    R(\delta')=r_{s_1-\delta'}\circ\ldots\circ r_{s_n-\delta'}.
  \]
  Observe that by definition $R_P=R_Q \circ R(
    \delta')\circ R(\delta)$. 

    Hence, we need to show that
  \[
    R(\delta')\circ R(\delta)= R'(\delta')\circ
  R'(\delta).
  \]
  At this level, $Q$ being a palindrome is not important, so we consider a general 
  \[
    Q=\eps^n_{e_n} \ldots
  \eps^n_1 \delta^n_{d_n} \ldots \delta^n_1
  \eps^{n-1}_{e_{n-1}}\ldots\ldots \eps^1_1 \delta^1_{d_1}\ldots\delta^1_1.
  \]
  Note that the number of odd reflections in $R(\delta')\circ
  R(\delta)$ is even. Hence, it makes sense to talk of the pairs of the first and the
  last reflection, the second and next to last, \etc 

  In particular, there is a central pair, which is connected by a string of even reflections of the form \eqref{eq:elementPair}.
  After removing this pair by Equation \eqref{eq:elementPair}, the next central pair is again only connected by even reflections, of
  the same form \eqref{eq:elementPair}. Iterating this procedure, we remove all odd
  reflections in pairs, and indeed find
  \[
    R(\delta')\circ R(\delta)= R'(\delta')\circ
  R'(\delta),
  \]
  as claimed. This completes the inductive proof of the first implication.

  For the converse, assume that there is a mismatching pair. That is, $P=A \delta B
  \eps \tilde A$ where $A$ and $\tilde A$ are of the same length, $A\tilde A$ is a (possibly empty) palindrome, and $B$ is any chain. Then $R_P = R_B
  \circ R(\eps) \circ R(\delta)\circ R(A,\tilde A)$, where after
  applying $R(X)$, 
  the $X$-components are in their final form, \ie will not be changed any
  more by consecutive reflections. So, by applying the first part of
  this proof, we have 
  \[
    R_P = R_B \circ R(\eps) \circ R(\delta)\circ R'(A,\tilde A).
  \]

  Denote by $\delta^*\in\ger a$ the dual element of $\delta$ in the $(\delta_i,\eps_j)$ basis. Then for any $\lambda \in \mf a^*$, we have 
  \[
    \delta^*(R_P \lambda) = \delta^*(R(\delta) \lambda).
  \]
  In particular, for 
  \[
    \lambda =-\sum_i \delta_i \in \mf a^*,
  \]
  where the sum runs over all $i$, all odd reflections in the chain $R(\lambda)$ act effectively on $\lambda$. Hence, we obtain 
  \[
    \delta^*(R(\delta) \lambda) < -1 = \delta^*(R_P' \lambda), 
  \]
  so that $R_P \neq R'_P$, proving the proposition. 
\end{PRF}

\begin{CRL}
  If $P$ is a palindrome, then $R_P$ acts as if the odd
  reflections were even, \ie in the same way as the element
  of the orthogonal group $O(\mf h^*)$ obtained from $R$ by replacing all odd reflections
  $r_\alpha$ by the ordinary even (linear) reflections perpendicular to $\alpha$.
\end{CRL}

The proof of \REF{pro}{mainSelfDual} also shows the following.

\begin{CRL}[prepostfix]
  In any $\delta\eps$-chain $AB\tilde A$, where $A$ and $\tilde A$ have equal length and $A\tilde A$ is a palindrome, one can drop the odd
  reflections from the palindrome part: 
  \[
  R_{AB\tilde A} = R_B\circ R(A,\tilde A) = R_B \circ R'(A,\tilde A).
  \]
\end{CRL}

\begin{CRL}[sd]
  All finite dimensional highest weight representations of $\gl^{p+q|r+s}$ with highest weight $\lambda$
  that vanishes on $\mf h \cap \mf k$, in particular all
  finite-dimensional spherical representations, are self-dual. 
\end{CRL}

\begin{PRF}
  Let $\lambda \in \mfa^*$ be a highest weight 
 and $R$ the chain of simple
  reflections mapping the highest to the lowest weight, defined in \REF{lmm}{BubbleSort}.
  By the special form of $\lambda$ in Equation \eqref{eq:lambda-astar}, only the
  palindrome part of $R$ is effective. Hence, we conclude that 
  \[
    R(\lambda)=R'(\lambda)=-\lambda
  \]
  by \REF{crl}{prepostfix} and Equation \eqref{eq:lambda-astar}.
\end{PRF}

\appendix
\section{Points and integrals}\label{app:point}

\subsection{The formalism of generalised points}

For any \emph{cs} manifold $X$, a point can be thought of as a morphism $* \to X$. However, if $X$ is not a manifold, then functions and morphisms on $X$ are not fully determined by their values on such ordinary points. To deal with this, the notion of points has to be extended.

Indeed, an \Define{$S$-valued point} (where $S$ is another \emph{cs} manifold) is defined to be a morphism $x:S\to X$. One may view this as a `parametrised' point. Suggestively, one writes $x\in_SX$, and denotes the set of all $x\in_SX$ by $X(S)$. For any morphism $f:X\to Y$, one may define a set-map $f_S:X(S)\to Y(S)$ by 
\[
  f_S(x)\defi f(x)\defi f\circ x\in_SY\mathfa x\in_SX.
\]
Clearly, the values $f(x)$ completely determine $f$, as can be seen by evaluating at the \Define{generic point} $x=\id_X\in_XX$.

In fact, more is true. The following statement is known as Yoneda's Lemma \cite{maclane}: Given a collection of set-maps $f_S:X(S)\to Y(S)$, there exists a morphism $f:X\to Y$ \scth $f_S(x)=f(x)$ \fa $x\in_SX$ if and only if 
\[
  f_T(x(t))=f_S(x)(t)\mathfa t:T\to S.
\]

The above facts are usually stated in the following more abstract form: For any $X$, we have a set-valued functor $X(-):\cat C^{op}\to\Sets$, where $\cat C$ is the category of \emph{cs} manifolds, and the set of natural transformations $X(-)\to Y(-)$ is naturally bijective to the set of morphisms $X\to Y$. Thus, the \Define{Yoneda embedding} $X\mapsto X(-)$ from $\cat C$ to the category $[\cat C^{op},\Sets]$ of set-valued functors, is a fully faithful functor.

The Yoneda embedding preserves products \cite{maclane}, so it induces a fully faithful embedding of the category of \emph{cs} Lie supergroups into the category $[\cat C^{op},\cat{Grp}]$ of group-valued functors. In other words, $X$ is a \emph{cs} Lie supergroup if and only if for any $S$, $X(S)$ admits a group law, which is compatible with base change in the sense defined above.

\subsection{Berezinian integrals}

Let $X$ and $S$ be \emph{cs} manifolds and $p:X\to S$ a morphism. Then $p$ or $X$ is called a \Define{\emph{cs} manifold over $S$}, written $X/S$, if on some open cover $U_\alpha$ of $X$ lying over an open cover $V_\alpha$ of $S$, we have commutative diagrams
\begin{center}
  \begin{tikzcd}
    U_\alpha\rar{}\dar[swap]{p} & S\times Y\dar{p_1}\\
    V_\alpha\rar{}&S
  \end{tikzcd}
\end{center}
where the rows are open embeddings. Such a diagram is called a \Define{local trivialisation} of $X/S$. Usually, we will consider only products, but the general language will be efficient nonetheless. There is an obvious notion of morphisms over $S$, which we denote $X/S\to Y/S$.

A system of (local) \Define{fibre coordinates} is given by the system $(x^a)=(x,\xi)$ of superfunctions on some trivialising open subspace $U\subseteq X$ obtained by pullback along a local trivialisation from a coordinate system in the fibre $Y$.

If $X/S$ is a \emph{cs} manifold over $S$, then the \Define{relative tangent sheaf} is defined by 
\[
  \sh T_{X/S}\defi\ShGDer[_{p_0^{-1}\sh O_S}]0{\sh O_X,\sh O_X},
\]
the sheaf of superderivations of $\sh O_X$ which are linear over $\sh O_S$. It is a basic fact that $\sh T_{X/S}$ is a locally free $\sh O_X$-module, with rank equal to the fibre dimension of $X/S$. 

Let $X$ be a \emph{cs} manifold over $S$, and $\smash{\Omega^1_{X/S}}$ be the module of relative $1$-forms, by definition dual to $\sh T_{X/S}$. Then we define the sheaf of \Define{relative Berezinians} $\sh Ber_{X/S}$ to be the Berezinian sheaf associated to the locally free $\sh O_X$-module $\smash{\Pi\Omega^1_{X/S}}$ obtained by parity reversal. Furthermore, the sheaf of \Define{relative Berezinian densities} $\Abs0{\sh Ber}_{X/S}$ is the twist by the relative orientation sheaf, \ie 
\[
  \Abs0{\sh Ber}_{X/S}\defi\sh Ber_{X/S}\otimes_\ints or_{X_0/S_0}.
\]

Given a system of local fibre coordinates $(x^a)=(x,\xi)$ on $U$, their coordinate derivations $\frac\partial{\partial x^a}$ form an $\sh O_X|_{U_0}$-module basis of $\smash{\sh T_{X/S}|_{U_0}}$, with dual basis $dx^a$ of $\smash{\Omega^1_{X/S}|_{U_0}}$. One may thus consider the distinguished basis 
\[
  \Abs0{D(x^a)}=\Abs0{D(x,\xi)}=dx_1\dots dx_p\frac{\partial^\Pi}{\partial\xi^1}\dots\frac{\partial^\Pi}{\partial\xi^q}
\]
of the module of Berezinian densities $\Abs0{\sh Ber}_{X/S}$, \cf \cite{Man88}.

If $X/S$ is a direct product $X=S\times Y$, then 
\[
  \Abs0{\sh Ber}_{X/S}=p_2^*\Parens1{\Abs0{\sh Ber}_Y}=\sh O_X\otimes_{p_{2,0}^{-1}\sh O_Y}p_{2,0}^{-1}\Abs0{\sh Ber}_Y.
\]
In particular, the usual Berezin integral over $Y$ of compactly supported Berezinian densities defines the integral over $X$ of a section of $(p_0)_!\Abs0{\sh Ber}_{X/S}$, where $(-)_!$ denotes the functor of direct image with compact supports \cite{iversen}. We denote the quantity thus obtained by 
\[
  \fibint[_S]{_X}\omega\in\Gamma(\sh O_S)\mathfa\omega\in\Gamma\Parens1{(p_0)_!\Abs0{\sh Ber}_{X/S}},
\]
and call this the \Define{fibre integral} of $\omega$.

We will, however, have to consider fibre integrals in a more general setting, beyond compact supports. Henceforth, we assume for simplicity that $X=S\times Y$. A \Define{fibre retraction} for $X$ is a morphism $r:Y\to Y_0$ which is left inverse to the canonical embedding $j:Y_0\to Y$, where $Y_0$ denotes the underlying manifold of $Y$. In case $S=*$, we just speak of a \Define{retraction} of $X$.

A special case, in which a distiguished retraction exists, is that of \emph{cs} Lie supergroups $G$. Indeed, $G$ is isomorphic to $G_0\times\ger g_\odd$, with the isomorphism given on $S$-valued points by 
\[
  (g,x)\mapsto g\exp_G(x).
\]
Since $p_1:G_0\times\ger g_\odd\to G_0$ is a retraction, we obtain one for $G$ by transport along this isomorphism. The retraction thus obtained is called the \Define{standard retraction} of $G$.

Return to $X/S$, with chosen retraction $r$. A system of fibre coordinates $(x,\xi)$ of $X/S$ is called \Define{adapted} to $r$ if $x=r^\sharp(x_0)$. Given an adapted system of fibre coordinates, we may write $\omega=\Abs0{D(x,\xi)}\,f$ and 
\[
  f=\sum_{I\subseteq\{1,\dotsc,q\}}(\id\times r)^\sharp(f_I)\,\xi^I
\]
for unique coefficients $f_I\in\Gamma(\sh O_{S\times Y_0})$, where $\dim Y=*|q$. Then one defines
\[
  \fibint[_{S\times Y_0}]{_X^r}\omega\defi\Abs0{dx_0}\,f_{\{1,\dotsc,q\}}\in\Gamma(\Abs0{\sh Ber}_{(S\times Y_0)/S}).
\]
Note that $\Abs0{\sh Ber}_{(S\times Y_0)/S}$ is $p_2^*$ of the sheaf of ordinary densities on the manifold $Y_0$, so we may write $\Abs0{dx_0}$.

This fibre integral only depends on $r$, and not on the choice of an adapted system of fibre coordinates. If the resulting relative density is absolutely integrable along the fibre $Y_0$, then we say that $\omega$ is \Define{absolutely integrable} with respect to $r$, and define
\[
  \fibint[_S]{_X^r}\omega\defi\fibint[_S]{_{S\times Y_0}}\Bracks3{\fibint[_{S\times Y_0}]{_X^r}\omega}\in\Gamma(\sh O_S).
\]
Both this quantity and its existence depend heavily on $r$. 

We shall use the language of $S$-valued points discussed above to manipulate inte\-grals of relative Berezinian densities in a hopefully more comprehensible formalism. This also gives a rigorous foundation for the super-integral notation common in the physics literature. 

If $f$ is a superfunction on $X=S\times Y$ and we are given some relative Berezinian density $\Abs0{Dy}$ on $X/S$, then we write
\[
  \int_Y\Abs0{Dy}\,f(s,y)\defi\fibint[_S]{_X^r}\Abs0{Dy}\,f.
\]
If the fibre retraction $r$ is understood, this is justified by the convention that the generic points of $S$ and $Y$ are denoted by $s$ and $y$, respectively. Moreover, it is easy to see that this notation behaves well under base change, since
\begin{equation}\label{eq:spec-int}
  \int_Y\Abs0{Dy}\,f(s(t),y)=\fibint[_T]{_X^r}(t\times\id)^\sharp(\Abs0{Dy}\,f)=t^\sharp\Bracks3{\fibint[_S]{_X^r}\Abs0{Dy}\,f}
\end{equation}
for any $t\in_TS$. This follows from the fact that the fibre retractions are respected by the morphism $s\times\id$.



\end{document}